\documentclass[12pt,leqno]{article}
\usepackage{amsfonts,amssymb,mathabx}
\usepackage{tikz}

\newcommand{\red}{}

\usetikzlibrary{matrix,arrows} 

 \def\medskipamount{12pt} \def\smallskipamount{6pt}

\newcounter{bitcount}
\newcommand{\bit}[1]{\addtocounter{bitcount}{1}\pagebreak[3]
\subsection{#1}\nopagebreak\setcounter{equation}{0}}

\renewcommand{\theequation}{\thesubsection .\arabic{equation}}
\renewcommand{\thesubsection}{\arabic{bitcount}}

\newcommand{\re}[1]{\mbox{\bf(\ref{#1})}}

\catcode`\@=\active
\catcode`\@=11
\renewcommand\@eqnnum{\hbox to .01pt{}\rlap{\bf \hskip -\displaywidth(\theequation)}}
\catcode`\@=12

\begin{document}


    \renewcommand{\topfraction}{0.9}	
    \renewcommand{\bottomfraction}{0.8}	
    \setcounter{topnumber}{2}
    \setcounter{bottomnumber}{2}
    \setcounter{totalnumber}{4}     
    \setcounter{dbltopnumber}{2}    
    \renewcommand{\dbltopfraction}{0.9}	
    \renewcommand{\textfraction}{0.07}	
    \renewcommand{\floatpagefraction}{0.7}	
    \renewcommand{\dblfloatpagefraction}{0.7}	


\catcode`\@=\active
\catcode`\@=11
\newcommand{\nc}{\newcommand}


\nc{\bs}[1]{ \addvspace{\medskipamount} \pagebreak[3]
\refstepcounter{equation}
\noindent {\bf (\theequation) #1.} \begin{em} \nopagebreak}

\nc{\es}{\end{em} \par \addvspace{\medskipamount} } 


\nc{\br}[1]{ \addvspace{\medskipamount} \pagebreak[3]
\refstepcounter{equation} 
\noindent {\bf (\theequation) #1.} \nopagebreak}


\nc{\er}{\par \addvspace{\medskipamount} }


\nc{\A}{\mathbb A}
\nc{\C}{\mathbb C}
\nc{\E}{\mathbb E}
\nc{\F}{\mathbb F}
\renewcommand{\H}{\mathbb H}
\renewcommand{\L}{\mathbb L}
\nc{\Q}{\mathbb Q}
\nc{\R}{\mathbb R}
\nc{\Z}{\mathbb Z}
\renewcommand{\P}{\mathbb P} 
\nc{\CC}{{\mathfrak C}}
\nc{\BB}{{\mathfrak B}}
\nc{\FF}{{\mathfrak F}}
\nc{\MM}{{\mathcal M}}
\nc{\NN}{{\mathcal N}}
\nc{\X}{{\mathcal X}}


\nc{\oper}[1]{\mathop{\mathchoice{\mbox{\rm #1}}{\mbox{\rm #1}}
{\mbox{\rm \scriptsize #1}}{\mbox{\rm \tiny #1}}}\nolimits}
\nc{\Ad}{\oper{Ad}}
\nc{\ad}{\oper{ad}}
\nc{\Aut}{\oper{H}}
\nc{\Auto}{\oper{Aut}}
\nc{\aut}{\oper{aut}}
\nc{\Bl}{\oper{Bl}}
\nc{\diag}{\oper{diag}}
\nc{\End}{\oper{End}}
\nc{\Env}{\oper{Env}}
\nc{\ev}{\oper{ev}}
\nc{\Ev}{\oper{Ev}}
\nc{\Ext}{\oper{Ext}}
\nc{\scrExt}{\oper{{\it Ext}}}
\nc{\Gr}{\oper{Gr}}
\nc{\Hom}{\oper{Hom}}
\nc{\scrHom}{\oper{{\it Hom}}}
\nc{\id}{\oper{id}}
\nc{\Isom}{\oper{Isom}}
\nc{\Tor}{\oper{Tor}}
\nc{\Pic}{{\oper{Pic }}}
\nc{\Proj}{{\oper{Proj }}}
\nc{\Spec}{{\oper{Spec }}}
\nc{\pr}{{\oper{pr}}}

\nc{\operlim}[1]{\mathop{\mathchoice{\mbox{\rm #1}}{\mbox{\rm #1}}
{\mbox{\rm \scriptsize #1}}{\mbox{\rm \tiny #1}}}}

\nc{\lp}{\raisebox{-.1ex}{\rm\large(}}
\nc{\rp}{\raisebox{-.1ex}{\rm\large)}}

\nc{\al}{\alpha}
\nc{\be}{\beta}
\nc{\la}{\lambda}
\nc{\La}{\Lambda}
\nc{\ep}{\varepsilon}
\nc{\ga}{\gamma}
\nc{\de}{\delta}
\nc{\si}{\sigma}
\nc{\om}{\omega}
\nc{\Om}{\Omega}
\nc{\Ga}{\Gamma}
\nc{\Si}{\Sigma}


\nc{\Left}[1]{\hbox{$\left#1\vbox to
    11.5pt{}\right.\nulldelimiterspace=0pt \mathsurround=0pt$}}
\nc{\Right}[1]{\hbox{$\left.\vbox to
    11.5pt{}\right#1\nulldelimiterspace=0pt \mathsurround=0pt$}}

\nc{\updown}{\hbox{$\left\updownarrow\vbox to
    10pt{}\right.\nulldelimiterspace=0pt \mathsurround=0pt$}}


\nc{\beqas}{\begin{eqnarray*}}
\nc{\actson}{\lefttorightarrow}
\nc{\co}{{\cal O}}
\nc{\cx}{{\C^{\times}}}
\nc{\down}{\Big\downarrow}
\nc{\Down}{\left\downarrow
    \rule{0em}{8.5ex}\right.}
\nc{\downarg}[1]{{\phantom{\scriptstyle #1}\Big\downarrow
    \raisebox{.4ex}{$\scriptstyle #1$}}}
\nc{\eeqas}{\end{eqnarray*}}
\nc{\fp}{\raisebox{.1ex}{     $\Box$} \par \addvspace{\smallskipamount}}
\nc{\G}{{\mathbb G}}
\nc{\lrow}{\longrightarrow}
\nc{\modmod}{/\! /}
\nc{\pair}{\! \cdot \!}
\nc{\pf}{\noindent {\em Proof}}
\nc{\sans}{\, \backslash \,}
\nc{\SG}{S_G} 
\nc{\st}{\, | \,}

\nc{\beq}{\begin{equation}}
\nc{\eeq}{\end{equation}}

\nc{\V}{{\rm V}}

\nc{\vv}[1]{{\bf ({\small V}.#1)}}

\def\vcentcolon{\mathrel{\mathop\ordinarycolon}}
\begingroup
  \catcode`\:=\active
  \lowercase{\endgroup
    \ifnum\mathcode`\:=32768\relax
      \let\ordinarycolon= :%
    \else
      \mathchardef\ordinarycolon\mathcode`\: %
    \fi
    \let :\vcentcolon}

\nc{\deq}{\vcentcolon =}

\renewcommand{\t}{{\mathfrak t}}
\renewcommand{\O}{{\mathcal O}}

\nc{\g}{{\mathfrak g}}
\nc{\m}{{\mathfrak m}}
\nc{\M}{{\mathfrak M}}

\hyphenation{para-met-riz-ing sub-bundle co-char-ac-ter non-trivial eval-u-a-tion para-metrized}

\catcode`\@=12

\noindent
{\LARGE \bf 
Compactifications of reductive groups}\vspace{7pt} \\
{\LARGE \bf as moduli stacks of bundles}
\\

\noindent {\bf Johan Martens } \\ 
School of Mathematics and Maxwell Institute, The
    University of Edinburgh\\ 
James Clerk Maxwell Building,
    Peter Guthrie Tait Road, Edinburgh EH9 3FD, Scotland\\
{\tt johan.martens@ed.ac.uk}

\smallskip

\noindent {\bf Michael Thaddeus } \\
Department of Mathematics, Columbia University\\ 
2990 Broadway, New York, NY 10027, USA\\
{\tt thaddeus@math.columbia.edu}

\renewcommand{\thefootnote}{}
\footnotetext{JM was supported by QGM (Centre for Quantum Geometry of Moduli
  Spaces) funded by the Danish National Research Foundation, and
  partially supported by FAPESP grant 2009/05136--2.}  
\footnotetext{MT was partially
  supported by NSF grants DMS--0401128 and DMS--0700419.}

\smallskip

\begin{quote}
{\small
{\bf Abstract.}  Let $G$ be a split reductive group.  We introduce the
moduli problem of {\it bundle chains\/} parametrizing framed principal
$G$-bundles on chains of lines.  Any fan supported in a Weyl chamber
determines a stability condition on bundle chains.  Its moduli stack
provides an equivariant toroidal compactification of $G$.  All toric
orbifolds may be thus obtained.  Moreover, we get a canonical
compactification of any semisimple $G$, which agrees with the
wonderful compactification in the adjoint case, but which in other
cases is an orbifold.  Finally, we describe the connections with Losev-Manin's
spaces of weighted pointed curves and with Kausz's 
compactification of $GL_n$.}
\end{quote}



\bigskip

\noindent
Let $G$ be a split reductive group over a field $k$.
Describing an equivariant compactification of $G$ with smooth
orbit closures is in some sense an old problem, with its roots
in the construction of complete quadrics and complete
collineations by Chasles, Schubert, Zeuthen, and others in the
19th century \cite{laksov}.  {\red Moduli problems
compactifying the classical groups were introduced by Semple
\cite{semple1, semple2, semple3}, Kleiman-Thorup \cite{kt}, and Kausz
\cite{kausz, kausz2} in the 20th century.}  In this paper we
embark on a new approach, {\red for general $G$}, using a
moduli problem {\red involving principal bundles on chains of
lines}.  This naturally leads us to construct toroidal
compactifications which, rather than being schemes, in general
are orbifolds, that is, smooth tame stacks.

The best solutions in the category of schemes are in the two extreme
cases where $G$ is either abelian or has trivial center.  The first is
provided by the theory of toric varieties.  The second appears in two
influential papers of De Concini-Procesi \cite{dp1,dp2}, where a
so-called {\it wonderful compactification} of a reductive group with
trivial center is introduced and the toroidal compactifications lying over
it are classified. 

The wonderful compactification has many attractive properties.  Most
notably, it has finitely many $G \times G$ orbits, each of whose
closures is smooth.  However, the approach of De Concini-Procesi
requires trivial center.  It is not immediately clear even what the
wonderful compactification of $SL_n$ should be.

An intriguing hint was dropped by Tonny Springer in his 2006 ICM
address \cite{springer}.  He pointed out that an arbitrary semisimple
$G$ may be compactified by taking the normalization in the function
field of $G$ of the wonderful compactification $\overline{G/Z_G}$.
This is known in the theory of spherical varieties as the canonical
embedding \cite{pezzini}.  Spherical varieties provide many
compactifications of reductive groups, but no canonical smooth
compactification of a semisimple group.

Indeed, the canonical embedding already fails to be smooth for $G=
SL_n$.  However, it has only finite quotient singularities and is
therefore the variety underlying an orbifold or smooth tame stack.
This suggests that it represents a moduli problem whose objects have
finite automorphism groups.
  
A clue as to the nature of this moduli problem came from a 1999 paper of
the second author \cite{ccr}. There it was shown that the wonderful
compactification $\overline{PGL_n}$ may be realized as the locus of
genus 0, degree $n$ stable maps to $Gr(n,2n+k)$ passing through two
general points.  Since this holds for any $k \geq 0$, we may
fancifully view this as a locus in the (nonexistent) space of genus 0,
degree $n$ stable maps to $BGL_n$.  Less fancifully, it still
parametrizes a family of genus 0 curves equipped with a vector bundle
framed at two points.  In fact, all of the curves in question turn out
to be {\it chains\/} of lines, so we might refer to the objects
parametrized as {\it framed bundle chains}.

This picture extends to the setting of principal $G$-bundles for any
split reductive $G$.  In fact, the structure group of any such bundle
on a chain reduces to a maximal torus $T$ \cite{kansas,teodorescu},
just as was shown by Grothendieck \cite{gk} and Harder \cite{harder} for the
projective line.  The isomorphism classes therefore turn out to be
parametrized by 1-parameter subgroups of $T$.  But when framings at
two marked points are introduced, continuous parameters appear.

A key insight was Charles Cadman's suggestion that everything be made
equivariant for an action of $\G_m$, the multiplicative group of the
field.  In retrospect, this seems natural by analogy with the relative
stable maps of Li \cite{li}, in which equivariant chains also loom
large.  Li was in turn inspired by the work of
Gieseker \cite{gieseker}, generalized by
Nagaraj-Seshadri \cite{nagaraj-seshadri} and Schmitt \cite{schmitt},
on vector bundles over semistable models of a nodal curve.  Their work
has some features in common with the present paper, but because of the
new element of $\G_m$-equivariance, the relationship is not
transparent.

\medskip

\noindent {\it The moduli problem.}  To be specific, we study the
moduli problem of {\it framed bundle chains}, defined to consist
of the following: (a) a chain of projective lines, of any length, with
the multiplicative group $\G_m$ acting with weight 1 on every
component; (b) a principal $G$-bundle on this chain, with a lifting of
the $\G_m$-action to it; (c) $\G_m$-invariant trivializations of this
bundle at the two smooth points $p_\pm$ fixed by $\G_m$.
The stack of such objects is neither separated
nor of finite type, but this may be cured by imposing something like a
stability condition, as follows.  

Thanks to the aforementioned reduction to $T$, those bundles on a
chain with $n$ nodes that admit a framing are classified (up to a Weyl
group action) by a sequence of 1-parameter subgroups $\be_1, \dots,
\be_n$ of $T$.  At the $i$th node, $\be_i$ indicates the weight of the
action on the fiber over the node.  We call $\be = (\be_1, \dots,
\be_n)$ the {\it splitting type}.

The stability condition depends on a choice of a simplicial fan $\Si$
supported on a Weyl chamber and a choice of lattice points $\be_1,
\dots, \be_N$, one on each ray of the fan.  A bundle is
defined in \re{m} to be $\Si$-{\it stable} if the components of its
splitting type span a cone in $\Si$ and form a subsequence of $\be_1,
\dots, \be_N$ in the correct order.  This bounds the length of the chain.

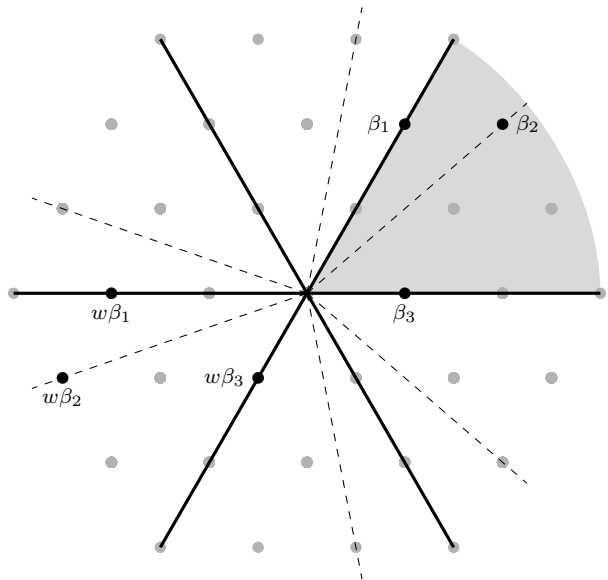
\begin{figure}[t]
\begin{center}
\begin{tikzpicture}

\newcommand{\dubbellinks}[1]
{\draw[line width=1pt](0,0)-- (.5,2.4); 
\draw[line width=1pt](.5,1.6)-- (0,4); 
\filldraw (.078,.4) circle (2pt);
\filldraw (.078,3.6) circle (2pt);
\draw (.15, .4) node[anchor=  east] {};
\draw (.15, .4) node[anchor=  west] {};
\draw (.15, 3.6) node[anchor=  east] {};
\draw (.15, 3.6) node[anchor=  west] {};
\draw ( .4,2) node[anchor = west]{#1};
}

\newcommand{\trippel}[2]
{
\draw[line width=1pt](0,0)--(-.5,1.6);
\draw[line width=1pt](-.5,1.2)--(0,2.8);
\draw[line width=1pt](0,2.4)--(-.5,4);
\filldraw (-0.13,.4) circle (2pt);
\filldraw (-0.37,3.6) circle (2pt);
\draw (-.1, .4) node[anchor=  east] {};
\draw (-.1, .4) node[anchor=  west] {};
\draw (-.3, 3.6) node[anchor=  east] {};
\draw (-.3, 3.6) node[anchor=  west] {};
\draw (-.4,1.466) node[anchor=west]{#2};
\draw (-.1,2.533) node[anchor=west]{#1};
}

\begin{scope}[xshift=-1.8cm, yshift=0cm, scale=1.3]

\fill[gray!30!white] (0,0) -- (3,0) .. controls (3,1) and (2.33, 2.1) .. (1.5, 2.6) -- cycle ;
\foreach \x in {-2,...,2}
\foreach \y in { -1.73, 0, 1.73}
 \foreach \z in {-2.5, -1.5,...,2.5}
\foreach \w in {-.866,.866}
{\filldraw[gray!60!white]  	(\x  , \y ) 		circle 	(1.5pt);
\filldraw[gray!60!white]   (\z , \w )	circle 	(1.5pt);}

\filldraw[gray!60!white]  	(.5  , 2.6 ) circle 	(1.5pt);
\filldraw[gray!60!white]	(-.5  , 2.6 ) circle 	(1.5pt);
\filldraw[gray!60!white]  	(-.5  , -2.6 ) circle 	(1.5pt);
\filldraw[gray!60!white]  	(.5  , -2.6 ) circle 	(1.5pt);
\filldraw[gray!60!white]  	(1.5  , 2.6 ) circle 	(1.5pt);
\filldraw[gray!60!white] 	(-1.5  , 2.6 ) circle 	(1.5pt);
\filldraw[gray!60!white]  	(1.5  , -2.6 ) circle 	(1.5pt);
\filldraw[gray!60!white]  	(-1.5  , -2.6 ) circle 	(1.5pt);
\filldraw[gray!60!white]  	(-3  , 0 ) circle 	(1.5pt);
\filldraw[gray!60!white]  	(3  , 0 ) circle 	(1.5pt);

\filldraw [black] 	(1,1.73 ) circle 	(1.5pt);
\filldraw [black] 	(2,1.73 ) circle 	(1.5pt);
\filldraw [black] 	(1,0 ) circle 	(1.5pt);
\filldraw [black] 	(-.5,-.866 ) circle 	(1.5pt);
\filldraw [black] 	(-2.5,-.866) circle 	(1.5pt);
\filldraw [black] 	(-2,0 ) circle 	(1.5pt);

\draw[very thick] (-3,0) -- (3,0);
\draw[very thick] (-1.5,-2.6) -- (1.5,2.6);
\draw[very thick] (1.5,-2.6) -- (-1.5,2.6);
\draw[dashed] (0,0) -- (2.25,1.946);
\draw[dashed] (0,0) -- (2.25,-1.946);
\draw[dashed] (0,0) -- (0.5625,2.925);
\draw[dashed] (0,0) -- (0.5625,-2.925);
\draw[dashed] (0,0) -- (-2.8125,.97425);
\draw[dashed] (0,0) -- (-2.8125,-.97425);

\draw (1,1.73) node[anchor = east] {$\scriptstyle{\beta_1}$};
\draw (2,1.73) node[anchor = west] {$\scriptstyle{\beta_2}$};
\draw (1,0) node[anchor = north] {$\scriptstyle{\beta_3}$};
\draw (-2,0) node[anchor = north] {$\scriptstyle{w\beta_1}$};
\draw (-2.5,-.866) node[anchor = north] {$\scriptstyle{w\beta_2}$};
\draw (-.5,-.866) node[anchor = east] {$\scriptstyle{w\beta_3}$};
\end{scope}

\begin{scope}[xshift=-8.5cm,yshift=-9cm, scale=.8]
\draw (-1,5) node[anchor = west] {\phantom{xxxx}\small $\Si$-stable bundle chains include};

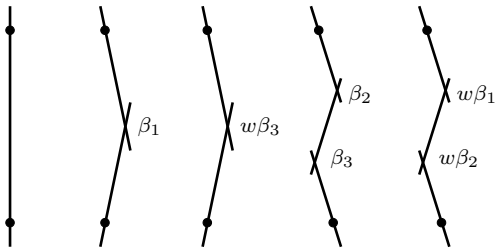
\begin{scope}[xshift=0cm, yshift= 2cm]
\draw[line width=1pt](0,2)--(0,-2);
\filldraw (0,-1.6) circle (2pt);
\draw (0.1, -1.6) node[anchor=  east] {};
\draw (0.1, -1.6) node[anchor=  west] {};
\filldraw (0,1.6) circle (2pt);
\draw (0.1, 1.6) node[anchor=  east] {};
\draw (0.1, 1.6) node[anchor=  west] {};
\end{scope}

\begin{scope}[xshift=1.5cm,yshift=0]
\dubbellinks{$\scriptstyle{\beta_1}$}
\end{scope}

\begin{scope}[xshift=3.2cm, yshift=0]
\dubbellinks{$\scriptstyle{w\beta_3}$}
\end{scope}

\begin{scope}[xshift=5.5cm,yshift=0]
\trippel{$\scriptstyle{\beta_2}$}{$\scriptstyle{\beta_3}$}
\end{scope}

\begin{scope}[xshift=7.3cm,yshift=0]
\trippel{$\scriptstyle{w\beta_1}$}{$\scriptstyle{w\beta_2}$}
\end{scope}

\end{scope}

\begin{scope}[xshift=.5cm,yshift=-9cm, scale=.8]
\draw (-2,5) node [anchor = west]{\phantom{xxxx}\small $\Si$-unstable bundle chains include};

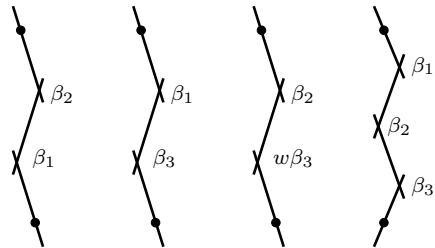
\begin{scope}[xshift=0cm,yshift=0]
\trippel{$\scriptstyle{\beta_2}$}{$\scriptstyle{\beta_1}$}
\end{scope}

\begin{scope}[xshift=2cm,yshift=0]
\trippel{$\scriptstyle{\beta_1}$}{$\scriptstyle{\beta_3}$}
\end{scope}

\begin{scope}[xshift=4cm,yshift=0]
\trippel{$\scriptstyle{\beta_2}$}{$\scriptstyle{w\beta_3}$}
\end{scope}

\begin{scope}[xshift=5.5cm, yshift= 0cm]
{\draw[line width=1pt](0,0)-- (.5,1.2); 
\draw[line width=1pt](.5,0.8)-- (0,2.2); 
\draw[line width=1pt](0,1.8)-- (.5,3.2); 
\draw[line width=1pt](.5,2.8)-- (0,4); 
\filldraw (.165,.4) circle (2pt);
\filldraw (.165,3.6) circle (2pt);
\draw (.22, .4) node[anchor=  east] {};
\draw (.22, .4) node[anchor=  west] {};
\draw (.22, 3.6) node[anchor=  east] {};
\draw (.22, 3.6) node[anchor=  west] {};}
\draw (.4, 3) node[anchor=  west] {$\scriptstyle{\beta_1}$};
\draw (0, 2) node[anchor=  west] {$\scriptstyle{\beta_2}$};
\draw (.4, 1) node[anchor=  west] {$\scriptstyle{\beta_3}$};
\end{scope}

\end{scope}

\end{tikzpicture}\\ \
\end{center}

\emph{\small A stacky fan $\Si$ supported in the positive Weyl
  chamber of $G=PGL_3$ (shaded in gray),
as well as the complete fan $W\Si$.  The action of a single element $w\in W$
is indicated.  Ordering the rays of $\Si$ determines a stability
condition on bundle chains, as illustrated by the examples.}
\begin{center}
\rule{5cm}{.1mm}
\end{center}

\end{figure}

The substack representing $\Si$-stable bundle chains is then proved to
be a smooth separated tame Artin stack of finite presentation, which
is proper if $\Si$ covers the entire Weyl chamber.  When the ground
field is $\mathbb C$, this means that it is a complex orbifold
(smooth, Hausdorff, and second countable), which is compact if $\Si$
covers the Weyl chamber.

When $G$ is a torus, this gives a modular interpretation to all toric
orbifolds.  On the other hand, when $G$ is semisimple, then taking $\Si$ to be the Weyl
chamber itself leads to a canonical choice of orbifold
compactification.  If $G$ moreover has trivial center, this is the familiar
wonderful compactification, but if not, it has nontrivial orbifold
structure except for $G=SL_2$.  
In general, it is a smooth stack whose coarse moduli space is an
arbitrary toroidal compactification, with finite quotient
singularities, of an arbitrary reductive group.

\medskip

\noindent {\it Summary of the paper.}  The paper
is organized as follows.  In \S1 we state the basic facts about chains
of projective lines and line bundles over them.  We give an explicit
atlas for the smooth Artin stack $\CC$ parametrizing chains {\red with
two marked points}, and we describe a canonical action {\red of} the
multiplicative group on any family of such chains.  In \S2 we
introduce {\it bundle chains}, our main objects of study.  These are
principal $G$-bundles on a chain of projective lines.  We classify
them up to isomorphism and study their deformation theory.  Then
in \S3 we prove that bundle chains are parametrized by a smooth Artin
stack $\BB$.  To get a tame or Deligne-Mumford stack $\MM$, we must
make three adjustments: rigidify to get rid of an overall
$\G_m$-action, add framings of our chains at two points, and impose a
mild open condition.  To get a separated stack $\MM(\Si)$ of finite
type (that is, an orbifold), one must go further and impose the
stability condition described above, determined by a stacky fan $\Si$.
This is carried out in \S4.

We then turn to an alternative construction of the moduli stack as a
global quotient of an algebraic monoid by a torus action.  This is
based on Vinberg's construction \cite{v} of the wonderful
compactification of an adjoint group in {\red the same} manner.  In \S5 we review
Vinberg's construction and show that it is a geometric invariant
theory quotient.  Then in \S6 we hybridize it with Cox's construction
\cite{cox} of toric varieties as global quotients of affine spaces.
As a result, $\MM(\Si)$ is expressed in \S7 as a geometric quotient
$S_{G,\be} \modmod \G_\be$ for a certain algebraic monoid $S_{G,\be}$
acted on by a torus $\G_\be$.  When it is projective, it is a
geometric invariant theory quotient.

The case where $\Si$ has a single cone covering the entire Weyl
chamber gives a canonical orbifold compactification of a semisimple
group.  In \S8 we show that its coarse moduli space agrees with the
space proposed by Springer.  

The final sections describe a few connections between our construction
and related ideas in the literature.  In \S9 we explain where the
coarse moduli spaces $M(\Si)$ of our stacks fall in the classification
of spherical varieties.  In \S10 we discuss the relationship with the
moduli problems represented by toric varieties and studied by
Losev-Manin \cite{lm}.  
Finally, in \S11, we show how the
compactification of the general linear group given by Kausz
\cite{kausz} fits into our picture.

\smallskip

\noindent {\it Acknowledgements.}  
We are happy to thank Dan Abramovich, Jarod Alper, J\o rgen Ellegaard
Andersen, Sergey Arkhipov, Michel Brion, Charles Cadman, Brian Conrad,
Johan de Jong, Edward Frenkel, Jochen Heinloth, Mathieu Huruguen, Lisa
Jeffrey, Ciprian Manolescu, Bart Van Steirteghem, Burt Totaro, Angelo
Vistoli, Jonathan Wise, Chris Woodward, {\red and the referee} for
their indispensable suggestions and advice.  Part of this work was
done during visits of JM to Columbia University and MT to QGM; we
gratefully acknowledge the hospitality and support of both
institutions.
\smallskip

\noindent {\it Notation.}  Angle brackets are used for the nonnegative
span of a set.  That is, the set of linear combinations of elements of
$S$ with nonnegative rational coefficients is denoted by $\langle
S\rangle$.  Also, the diagonal in $S \times S$ is denoted by
$\Delta_S$.  Notation related to the reductive group $G$ and its
characters and cocharacters is summarized in \re{dd}.  The center of
$G$ is denoted $Z_G$.  For a cocharacter $\la$, we denote by $t^\la$
the value of the corresponding 1-parameter subgroup on $t \in \G_m$.
This has the virtue of yielding $t^{\la + \la'} = t^\la t^{\la'}$.

\bit{Chains and the stack of chains}

Let $k$ be a field.  All schemes and stacks throughout will be over
$k$.

The objects of our moduli problem will be principal bundles with
reductive structure group over a chain of projective lines,
equivariant for the multiplicative group action, and trivialized at
the two endpoints.  For brevity, we will refer to them as {\em framed bundle
  chains}.

\br{Definition} The {\em standard chain with $n$ nodes\/} over a
field extension $K/k$ is the nodal curve $C_n$ each of whose $n+1$ irreducible
components is isomorphic to the projective line $\P^1_K$, with nodes
where $[0,1]$ in the $i-1$th line is glued to $[1,0]$ in the $i$th
line.  Its {\em endpoints} are $p_+ \deq [1,0]$ in the first line and
$p_- \deq [0,1]$ in the last line.  They are fixed points of
the {\em standard action} $\G_m \actson C_n$ given by $t \cdot
[u,v] \deq [u, tv]$ on each component. \er
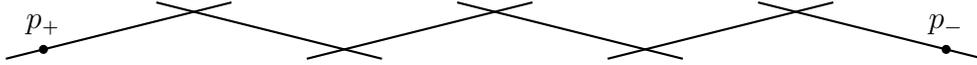
\begin{figure}[h]                                                               
\begin{center}                                                                  
\begin{tikzpicture}[scale=1]                                                    
\draw[thick] (1.5, .375) -- (-1.5, -.375);                                
\draw[thick] (.5, .375) -- (3.5, -.375);                                  
\draw[thick] (5.5, .375) -- (2.5, -.375);                                 
\draw[thick] (4.5, .375) -- (7.5, -.375);                                 
\draw[thick] (9.5, .375) -- (6.5, -.375);                                 
\draw[thick] (8.5, .375) -- (11.5, -.375);                                
\filldraw (-1, -.25) circle (1.5pt);                                            
\draw (-1, -.2) node[anchor= south] {$p_+$};                                    
\filldraw (11, -.25) circle (1.5pt);                                            
\draw (11, -.2) node[anchor= south] {$p_-$};                                    
\end{tikzpicture}                                                               
\end{center}                                                                    
\end{figure}      
\br{Remark} \label{c} 
We will be concerned with $\G_m$-equivariant principal bundles over $C_n$, or
equivalently, principal bundles over the stack $[C_n/\G_m]$.
Regarding line bundles, 
the following are easily verified:

(a) The map $\Pic [C_n/\G_m] \to \Z^{n+2}$ taking an equivariant
line bundle to its weights on the fixed points is an isomorphism.

(b) Denote by $\co(b_0
\st b_1 \st \dots \st b_{n+1})$ the line bundle corresponding to
$(b_0, \dots, b_{n+1}) \in \Z^{n+2}$ under this isomorphism.  
Then the dualizing sheaf is $\om_{C_n} \cong \co(-1\st 0 \st \dots \st
0 \st 1)$.  That is, $\om_{C_n} \cong \O(-p)$, where $p
\deq p_+ + p_-$.

(c) For $1 \leq i \leq j \leq n$, there is a section of $\co(b_0 \st \dots \st
b_{n+1})$ over $[C_n/\G_m]$ whose support is the $i$th through $j$th
components if and only if $b_{i-1} > 0 = b_i = \cdots = b_{j-1} >
b_j$, except that $b_{i-1}=0$ is allowed if $i=1$ and $b_j=0$ is allowed
if $j=n+1$: in which cases, the section is nonvanishing at $p_+$
and $p_-$, respectively.

\er

\br{Definition} 
\label{bb}
A {\em chain\/} over a scheme $X$ is an algebraic space
$C$, flat, proper,
and finitely presented over $X$, equipped
with two sections $p_\pm: X \to C$, such that every geometric fiber of
$C \to X$ is isomorphic to a standard chain.  \er

Since we choose no polarization on our chains, a flat
family of chains need not be projective, or even a scheme, just
an algebraic space.

\br{Example} 
\label{d} 
The {\it versal chain} $\C_n \to \A^n$ is defined recursively for any
$n \geq 0$ by blowing up as follows.  Let $\C_0 \deq \P^1$ with $p_+
\deq [1,0]$ and $p_- \deq [0,1]$.  Then, given $\C_n$, let $\C_{n+1}
\deq \Bl(\C_n \times \A^1, p_-(\A^n) \times 0)$ with $p_\pm$ defined
as the proper transforms of their counterparts on $\C_n$.  The
following properties are now easily verified:

(a) The obvious action $\G_m^{n+1} \actson \P^1 \times \A^n$ induces
an action $\G_m^{n+1} \actson \C_n$ in which the action of $\G_m
\times 1^n$ is the standard action on each fiber, while the action of
$1 \times \G_m^n$ lies over the obvious action on $\A^n$.


(b) Let $E_i$ be the proper transform in $\C_n$ of the exceptional
divisor of the $i$th blow-up (the one giving rise to $\C_i$), and let
$\L_i \deq \O(-E_i)$ over $\C_n$.  Also let $\L_0 \deq \O(p_+)$ and
$\L_{n+1} \deq \O(-p_-)$.  Then (i) $\G_m^{n+1}$ acts naturally
on $\L_i$, that is, $\L_i \in \Pic[\C_n / \G_m^{n+1}]$; (ii) the
restriction of $\L_i$ to the central fiber $[C_n/\G_m]$ is
$\O(0|\cdots|0|1|0|\cdots|0)$, where $1$ appears at the $i$th position;
(iii) the $j$th factor of $\G_m^n$ acts on the fiber of $\L_i$ over
$p_+(0)$ trivially and on that over $p_-(0)$ with weight
$-\delta_{i,j}$; (iv) the $\L_i$ freely generate the kernel of the
restriction $\Pic [\C_n/\G_m^{n+1}] \to \Pic[p_+(0)/\G_m^n]$.

(c) For any $S \subset \{1, \dots, n\}$, let $U_S \deq \{ x \in
\A^n \st x_i \neq 0 \mbox{ if } i \not\in S\}$.  Then the morphism
$\pi_S : U_S \to \A^{|S|}$ given by projection on the coordinates
in $S$ is a geometric quotient by a subtorus of $\G_m^n$,
inducing a $\G_m^n$-equivariant isomorphism of chains $\C_n |_{U_S}
\cong \pi_S^* \C_{|S|}$ under which the line bundles
$\L_i$ correspond.
\er

The following statement and proof are well known (indeed, the stack
{\red described} is an open substack of the stack ${\mathfrak M}_{0,2}$
of prestable curves {\red with two smooth marked points}), but we
include them for completeness.

\bs{Theorem} The category of chains (and isomorphisms thereof) is a smooth
Artin stack $\CC$ with atlas $\bigsqcup_{n \geq 0} \C_n:
\bigsqcup_{n \geq 0} \A^n \to \CC$. 
\label{h}
\es


\pf.  We must prove three things: (1) that $\CC$ is a stack,
namely (a) a category fibered in groupoids, such that (b) \'etale
descent data are effective and (c) automorphisms are a sheaf; (2) that
the diagonal $\CC \to \CC \times \CC$ is representable, separated, and
finitely presented; and (3) that the atlas stated is surjective and
smooth, so that the stack is algebraic.  The smoothness of $\CC$ then
follows directly from that of $\A^n$ \cite[4.7.1]{lmb}, {\red as local
properties by definition hold for an algebraic stack if and only if
they hold for an atlas.}

Step (1).  That $\CC$ is (a) a category fibered in groupoids is
immediate from the existence of fibered products in the category of
algebraic spaces \cite[02X2]{stacks}, \cite[II 1.5]{knutson}.  To prove
(b), observe that since a chain is of finite type over its base,
\'etale descent data are effective for chains \cite[04UD]{stacks}, and
the property that the geometric fibers are standard chains is
preserved by \'etale descent, since for algebraically closed $K$,
\'etale $X' \to X$, and $\Spec K \to X$ a geometric point, $X'
\times_X \Spec K$ is \'etale over $\Spec K$, hence a disjoint union of
copies of $\Spec K$, so there is a section $\Spec K \to X' \times_X
\Spec K$, identifying each geometric fiber of a descended chain $C \to
X$ with some fiber of the original $C' \to X'$.  That (c) isomorphisms
between chains $C_1, C_2 \to X$ constitute a sheaf also follows from
descent for algebraic spaces, since any automorphism can be identified
with its graph, a closed algebraic subspace of $C_1 \times_X C_2$.
Hence $\CC$ is a stack.

Step (2).  The requisite properties of the diagonal are verified
exactly as by Fulghesu \cite[1.9]{fulghesu} and Hall \cite[\S3]{hall}.

Step (3).  Finally, we must show that $\bigsqcup_{n \geq 0} \A^n \to
\CC$ is smooth and surjective on geometric points.  The surjectivity
is obvious, since the standard chain with $n$ nodes appears over the
origin of $\A^n$.  The smoothness is by definition \cite[3.10.1]{lmb}
a matter of showing, for any chain $C \to X$ inducing $X \to \CC$,
that $\A^n \times_\CC X \to X$ is smooth.  By noetherian reduction we
may assume $X$ locally noetherian.  By passing to an \'etale cover and
using descent \cite[1.15]{vistoli}, \cite[II 3.{\red 2}]{knutson}, we may
also assume that $C$ is a scheme. Then it suffices to show \cite[IV
17.14.2]{ega} that it is locally of finite presentation and formally
smooth in the sense that a lift always exists in the following
diagram, where $A$ is an Artinian $K$-algebra for some field extension
$K/k$ and $J \subset A$ is an ideal with $A/J \cong K$:
\begin{center}
\begin{tikzpicture}[>=angle 60] 
\matrix (m) [matrix of math nodes, row sep=2em, 
column sep=1em, text height=1.5ex, text depth=0.25ex] 
{ \Spec A/J & &  \A^n \times_\CC X \\ 
\Spec A & & X.\\ }; 
\path[->] 
(m-1-1) edge (m-1-3);
\path[->] 
(m-1-1) edge (m-2-1); 
\path[->] 
(m-1-3) edge (m-2-3); 
\path[->]
(m-2-1) edge (m-2-3);
\path[dashed,->]
(m-2-1) edge (m-1-3);
\end{tikzpicture}
\end{center}

It is locally of finite presentation by the corresponding property of
the diagonal established above: since $\A^n \times_\CC X
\cong (\A^n \times X) \times_{\CC \times \CC}
\CC$, it is locally of finite presentation over $\A^n \times
X$ and hence over $X$.  

To show the formal smoothness, observe that a
lift consists of a morphism $\Spec A \to \A^n$ extending the given
$\Spec A/J \to \A^n$ and an isomorphism of the two chains
on $\Spec A$ pulled back from $\A^n$ and from $X$, extending the given
one on $\Spec A/J$.  But we know the following:

(a) The set of isomorphism classes of chains over $\Spec A$ extending a given
chain over $\Spec A/J$ is nonempty and is acted on transitively by 
$\Ext^1(\Omega, \O(-p))$, 
the space of first-order deformations of
$C_n$ with marked points $p = p_+ + p_-$.  See Sernesi \cite[2.3.4,
2.4.1, 2.4.8]{sernesi}.   The addition of the two marked points does
not affect the proofs.

(b) The set of morphisms $\Spec A \to \A^n$ extending the given $\Spec A/J
\to \A^n$ is acted on transitively (and freely) by the tangent space
to $\A^n$ at the image of $\Spec A/{\mathfrak m}$.  This is elementary.

(c) The natural map from the latter to the former is given, in terms
of these actions, by the derivative or Kodaira-Spencer map of $\C_n$.
Indeed, by Schlessinger's condition $H_2$ (which Sernesi calls $H_\ep$),
chains over $\Spec A \otimes_K K[\ep]/(\ep^2)$ naturally correspond to
pairs of chains over $\Spec A$ and $\Spec K[\ep]/(\ep^2)$ with an
isomorphism of the central fiber.  The actions of (a) and (b) are then
given by base change by the morphism $\Spec A \to \Spec A \otimes_K
K[\ep]/(\ep^2)$ taking $\ep$ to a generator of $J$.  This clearly
commutes with the map (b) $\to$ (a) induced by $\C_n$; but the map
induced on $\Spec K[\ep]/(\ep^2)$ is exactly the Kodaira-Spencer map.

Finally, the family $\C_n \to \A^n$ has surjective Kodaira-Spencer map at
every point.  Indeed, the space of first-order deformations is $\Ext^1(\Omega, \O(-p))
= H^0(\scrExt^1(\Omega, \O(-p)))$ as in Sernesi
\cite[1.1.11]{sernesi}.  This
has a basis
consisting of elements smoothing a single node, and moving along the
coordinate axes in $\A^n$ carries out these deformations.  \fp

\bs{Theorem}  
\label{i}
Let $C \to X$ be a chain.  There exists an unique $\G_m$-action on $C$
lying over the trivial action on $X$ and restricting to the standard
action on each geometric fiber.  Moreover, for any line bundle $L \to
C$, there is an unique lifting of this action to $L$ that acts
trivially on $p_+^*L$.  
\es

This canonical $\G_m$-action will be called the {\it uniform action}.


\pf.  Let the {\it standard equivariant chain with $n$ nodes} be the
standard chain equipped with the standard $\G_m$-action, and let an
{\it equivariant chain\/} be a chain equipped with a $\G_m$-action
lying over the trivial action on the base and restricting to the
standard action on each geometric fiber.  Then the category of
equivariant chains is again an Artin stack $\tilde \CC$.  Indeed, most
of the proof runs exactly parallel to that of \re{h}.  To complete
Step (2), observe that an equivariant isomorphism of chains is just an
ordinary isomorphism that satisfies a closed condition, namely that it
intertwine the two $\G_m$-actions.  Hence for any two equivariant
chains $C_1, C_2 \to X$, the corresponding morphism $X
\times_{\tilde\CC \times \tilde\CC} \tilde\CC \to X \times_{\CC \times
  \CC} \CC$ is a closed immersion.  Hence it is separated \cite[I
5.5.1]{ega} and of finite type \cite[I 6.3.4]{ega}.  By noetherian
reduction we may assume $X$ is locally noetherian.  {\red Since
$\CC \to \CC \times \CC$ is of finite type by Step (2) of \re{h},
$X \times_{\CC \times \CC} \CC \to X$ is of finite type \cite[I 6.3.4]{ega}.}  Then $X
\times_{\CC \times \CC} \CC$ is also locally noetherian \cite[I
6.3.7]{ega}, and hence our closed immersion is finitely presented
{\red \cite[IV 1.6.1]{ega}}.  Therefore the composition
$X \times_{\tilde\CC \times \tilde\CC}
\tilde\CC \to X$ is separated and finitely presented, so the same is
true of the diagonal $\tilde\CC \to \tilde\CC \times \tilde\CC$.

There is a forgetful morphism $\tilde\CC \to \CC$, and clearly it
suffices to show it is an isomorphism.  It is finitely presented,
since the same is true of $\tilde\CC \to \Spec k$ and of the diagonal
of $\CC \to \Spec k$.  It is universally injective and surjective
simply because every standard chain over every field admits an unique
$\G_m$-action isomorphic to the standard action.  Thanks to the
fundamental property of \'etale morphisms \cite[IV 17.9.1]{ega}, it
remains only to show that it is \'etale.  Since it is finitely
presented, it suffices to show that it is formally \'etale for
morphisms from Artinian $K$-algebras.  As in the proof of \re{h}, it
suffices to show that the derivative or Kodaira-Spencer map is an
isomorphism.  The first-order deformations are now
$\Ext^1_{[C_n/\G_m]}(\Om,\O(-p))$.  This may be shown to equal
$\Ext^1_{C_n}(\Om,\O(-p))^{\G_m}$ via the spectral sequence of the
hypercovering associated to the presentation $\G_m \times C_n
\rightrightarrows C_n$.  But $\Ext^1_{C_n}(\Om,\O(-p)) =
H^0(\scrExt^1_{C_n}(\Om,\O(-p)))$, and $\scrExt^1_{C_n}(\Om,\O(-p))$
is a sum of skyscraper sheaves supported at the nodes and acted on
trivially by $\G_m$.  Hence $\Ext^1_{C_n}(\Om,\O(-p))^{\G_m} =
\Ext^1_{C_n}(\Om,\O(-p))$ as desired.
  
To prove the last statement, let $\phi: \G_m \times C \to \G_m \times
C$ be given by $\phi(z,c) = (z,z\cdot c)$, acting on the second factor
as above.  To lift the action to $L$ means to find an isomorphism
$\phi: L \to L$ satisfying the obvious commutative diagram.  From the
classification of line bundles on standard chains \re{c}(a), it is
clear that there exists an action over each $\G_m \times \{ s \}$.  So
$L^{-1} \otimes \phi^*L$ is trivial on the fibers of the projection
$\id \times \pi: \G_m \times C \to \G_m \times X$, and hence the
direct image $\pi_*(L^{-1} \otimes \phi^*L)$ is a line bundle on $\G_m
\times X$.  It is canonically trivialized on $1 \times X$, hence
trivial.  So there is an isomorphism $L \cong \phi^*L$ over $\G_m
\times C$ restricting to the identity on $1 \times C$.  In fact any
such isomorphism satisfies the commutative diagram: after all, over
each $\G_m \times \{ s \}$, one isomorphism does, and any two such
isomorphisms differ by multiplication by a morphism $\G_m \to \G_m$,
which is of the form $z \mapsto z_0 z^n$ for some $z_0 \in \G_m$ and
$n \in \Z$.  Because the restriction to $1 \times C$ is
required to be the identity, $z_0 = 1$, so all such isomorphisms
differ by tensoring by a character of $\G_m$ and hence all satisfy the
diagram.  Tensoring by a further such character, one may arrange that
the action on $p_+^*L$ is trivial.~\fp

\vspace{15pt}

\bit{Bundle chains}

We now introduce our main objects of study.

\br{Notation}
\label{dd}
Let $G$ be a split reductive algebraic group over $k$, $T$ a split maximal
torus, $W \deq N(T)/T$ the Weyl group, $\V \deq \Hom(T,\G_m)$ the
character lattice, $\La \deq \Hom(\G_m,T)$ the cocharacter lattice,
$B$ the Borel subgroup corresponding to a choice of simple roots
$\al_i$, and $\La_\Q^+ \subset \La \otimes \Q$ the positive Weyl
chamber.  \er

\br{Definition}
A {\em principal $G$-bundle} $E$ over a scheme (or stack) $X$ is a
$G$-torsor locally trivial in the \'etale topology.  Its {\it adjoint
  bundle} $\Ad E$ is the associated bundle with fiber $G$, but with transition
functions given by conjugation (instead of multiplication) by the
transition functions of $E$.  Likewise, $\ad E$ is the associated vector bundle
with fiber ${\mathfrak g}$.
\er

\br{Definition}
\label{a}
A {\em $G$-bundle chain\/} over a scheme (or stack) $X$ is a chain $C
\to X$ equipped with a principal $G$-bundle $E$ over the stack
$[C/\G_m]$, where the action of $\G_m$ on $C$ is the uniform one of
\re{i}.  (Equivalently, this is a $\G_m$-equivariant principal
$G$-bundle over $C$.)  Likewise, a {\em framed bundle chain\/} is a
bundle chain $E$ equipped with ($\G_m$-invariant) trivializations of
$p_\pm^* E$.  \er

\br{Remark}
We shall be interested only in bundle chains where the restriction of
the bundle to each geometric fiber is {\em rationally trivial}, that
is, trivial on the generic point of each component.  This notion will
play a minor role, as it is automatic if $k$ is
algebraically closed of characteristic zero
\cite[1.9]{steinberg}, if $G$ is a torus 
\cite[III 4.9]{milne}, or if the bundle chain is framed \cite[1.1]{bn}.  
It is also clearly an open condition, as the deformation space
$H^1(\Spec K(t), \ad E)$ of a trivial bundle $E$ over a generic point
$\Spec K(t)$ of a chain certainly vanishes.  So it will not affect our
deformation theory arguments. \er

\br{Example} 
\label{f}
Suppose the chain is simply $C = \P^1 \times X \to X$.  Then any
framed bundle chain over $C$ must be trivial as a bundle, as we will
see in \re{b} below.  The framing at $p_+$ fixes a global trivialization,
from which the framing at $p_-$ differs by a morphism $X \to G$.
Hence the moduli space of such framed bundle chains is nothing
but $G$.  The moduli stacks we eventually study will compactify this
locus. \er

\br{Example} 
\label{ee}
Suppose the chain is a standard chain $C_n \to \Spec K$.
A $T$-bundle over $C_n$ is essentially an $r$-tuple of line bundles
and as such, according to \re{c}(a), is determined up to isomorphism
by the weights of the $\G_m$-action at the fixed points of $C_n$.  Let
$\be = (\be_1\dots, \be_n) \in \La^n$, let $F(\be)$ be the $T$-bundle whose weights at $p_\pm$ are 0 and whose weights at the nodes are $\be_1,\dots,\be_n$, and let
$E(\be)$ be the associated $G$-bundle.  Then
$E(\be)$ defines a bundle chain over $\Spec K$. \er

\br{Example} 
\label{g}
Again let $\be = (\be_1, \dots, \be_n) \in \La^n$.  Then {\it versal bundle
  chains} $\E(\be)$ and $\F(\be)$ over $\A^n$, with structure groups $G$ and $T$
respectively, are defined as follows. Let $\C_n \to \A^n$ be the
versal chain of \re{d}, $\L_1, \dots, \L_n \in \Pic [\C_n
/\G_m^{n+1}]$ as described there.  Then $\L_1^\times \times \cdots
\times \L_n^\times$ is a principal $\G_m^n$-bundle.  Regarding $\be 
\in \La^n$ as an element of $\Hom (\G_m^n, T)$, let $\F(\be)$ be
the associated $T$-bundle and $\E(\be)$ the associated $G$-bundle.
Then $\F(\be)$ and $\E(\be)$ are bundle chains over $[\A^n/\G_m^n]$.
\er

Regarded as bundle chains over $\A^n$, both $\F(\be)$ and $\E(\be)$
may of course be framed at $p_\pm$.  It is convenient to single out a
choice of {\it conventional framings} $\Psi_\pm: \A^n \times G \to
p_\pm^* \E(\be)$ as follows.  Since $\L_i = \O(-E_i)$, the
functions $1$ and $x_i$ constitute nowhere vanishing sections of
$p_+^*\L_i$ and $p_-^*\L_i$ respectively.  They induce framings of
$\L_1^\times \times \cdots \times \L_n^\times$ and hence of $\F(\be)$ and
$\E(\be)$.

Over the open set $x_1 \cdots x_n \neq 0$, all of the above bundles are
canonically trivial.  Relative to this trivialization, the
conventional framings may be expressed as $1$ and $x^\be \deq \prod
x_i^{\be_i}$ respectively.  Therefore, the conventional framings
are {\em not\/} both $\G_m^n$-equivariant: $\Psi_+$ is, but $\Psi_-$ is not.

\br{Definition} The {\it automorphism group\/} of $E$ is $\Ga(\Ad E)$,
cf.\ Brion \cite{brion2}.  However, the notation $\Auto E$ is reserved for
the {\it framed automorphism group}, that is, the subgroup of $\Ga
(\Ad E)$ preserving the framings at $p_\pm.$ Likewise, $\Auto C$
denotes the framed automorphisms of the chain $C$, namely those fixing
$p_+$ and $p_-$; and $\Auto(C,E)$ denotes the framed automorphisms of
the bundle chain, so that
$$1 \lrow \Auto E \lrow \Auto (C,E) \lrow \Auto C.$$
\er

\br{Definition} For a subgroup $\iota :H \to G$, a {\em reduction of
  structure group} of a principal $G$-bundle $E \to X$ is a principal
$H$-bundle $F\to X$ and an isomorphism $F_\iota \cong E$, where
$F_\iota$ is the associated $G$-bundle $(F \times G)/H$.  Two
reductions are {\it isomorphic} if there is an isomorphism of the two
$H$-bundles making the obvious triangle commute.  Isomorphism classes
of reductions then correspond to sections of the associated
$G/H$-bundle $E/H \to X$.  For example, when $G = GL_r$ and the
maximal torus $T$ consists of diagonal matrices, a reduction to $T$ is
equivalent to a splitting as a sum of line bundles.  \er

\bs{Theorem} 
Every rationally trivial $G$-bundle chain over $\Spec K$ admits a
reduction of structure group to $T$.  Its isomorphism class is unique
modulo the actions of $\Ga(\Ad E)$ and the Weyl group $W$.  \es

\pf. See another paper by the authors \cite[6.4]{kansas}. \fp

\bs{Corollary} 
\label{b}
Let $\bar{H}^1(X, G)$ denote the set of isomorphism classes of rationally
trivial principal $G$-bundles on $X$.  Then there is a natural
bijection $\bar{H}^1([C_n/\G_m],G) = \La^{n+2}/W$.  The $G$-bundles
admitting a framing correspond to $n+2$-tuples whose first and last
coordinates vanish, hence are in bijection with $\La^n/W$.  That is,
they are the bundles $E(\be)$ defined in {\rm \re{ee}}.
\es

\pf. By \re{c}(a), $\Pic [C_n/\G_m] \cong \Z^{n+2}$, and by Theorem
90 \cite[III 4.9]{milne} all line bundles are rationally trivial.
Hence $\bar{H}^1([C_n/\G_m]; T) \cong \La^{n+2}$.  The first statement
then follows from the theorem.  A bundle corresponding to
$(\be_0, \dots, \be_{n+1}) \in \La^{n+2}$ 
admits a framing if and only if the action of $\G_m$ at the endpoints
is trivial, that is, $\be_0 = \be_{n+1} = 0$. \fp

Let $\Om$ be the K\"ahler differentials on $C_n$, and let ${\mathcal
S}$ denote the sheaf on $C_n$ of $G$-invariant K\"ahler
differentials on the total space of a framed bundle chain $E$.  There
is then a short exact sequence of sheaves on $[C_n/\G_m]$
$$0 \lrow \Om \lrow {\mathcal S} \lrow \ad E \lrow 0,$$ in which the
three nonzero terms control the deformation theory of the chain $C_n$,
of the bundle chain $(C_n,E)$, and of the bundle $E$ respectively.
More precisely, let $T^0_{C_n}$, $T^1_{C_n}$, $T^2_{C_n}$ denote
first-order endomorphisms, deformations, and obstructions,
respectively, of the chain $[C_n/\G_m]$ {\red together with its two
marked points}.  (For simplicity, this notation suppresses any mention
of the {\red marked points} or of the $\G_m$-equivariance.)  Likewise, let
$T^i_{C_n,E}$ denote the corresponding spaces for the framed bundle
chain $E$ over $[C_n/\G_m]$, and let $T^i_E$ denote the corresponding
spaces for the framed bundle $E$ over the fixed base $[C_n/\G_m]$.

\bs{Proposition}
\label{q}
For $i=0,1,2$, there are natural isomorphisms
$$\, T^i_{C_n} = \Ext^i_{C_n}(\Om,\co(-p))^{\G_m}, \quad
\, T^i_{C_n,E} = \Ext^i_{C_n}({\mathcal S},\co(-p))^{\G_m}, \quad
\, T^i_{E} = \Ext^i_{C_n}(\ad E,\co(-p))^{\G_m}.$$
In particular all obstructions vanish.
\es

\pf. The proofs that $T^i_{C_n} =
\Ext^i_{[C_n/\G_m]}(\Om,\co(-p))$, and similarly for the other two
cases, are straightforward following Sernesi \cite[2.4.1,
3.3.11]{sernesi}.  But the latter agrees with the
$\G_m$-invariant part of $\Ext^i_{C_n}$.  This can be shown by
computing $\Ext^i_{[C_n/\G_m]}$ via the spectral sequence associated to
  the hypercovering of the presentation $\G_m \times C_n \rightrightarrows C_n$. 
\fp

Any bundle chain possesses a trivial 1-parameter group of
automorphisms, namely the $\G_m$-action itself.
Hence $\dim T^0_{C_n,E} \geq 1$.  The following result explains
when equality holds, that is, when there are no other infinitesimal
automorphisms.

\bs{Theorem} 
\label{j}
Let $\be = (\be_1, \dots, \be_n)$ and $E =
E(\be)$.  Then

{\rm (a)} $\dim T^0_E = 0$ if and only if $\be_1, \dots, \be_n$ lie in a
common Weyl chamber;

{\rm (b)} $\dim T^0_{C_n,E} = 1$ if and only if
$\be_1, \dots, \be_n$ lie in a
common Weyl chamber and are linearly independent in $\La \otimes k$;

{\rm (c)} $\dim \Auto (C_n, E) = 1$ if and only if
$\be_1, \dots, \be_n$ lie in a
common Weyl chamber and are linearly independent in $\La$.

\es

Conditions (b) and (c) are, of course, equivalent in
characteristic zero.  In positive characteristic the $\be_i$ may be
dependent in $\La \otimes k$ but not in $\La$.  In this case
$\Auto(C_n,E)$ will not be smooth.

\smallskip


\pf.  
(a): We have $T^0_E = \Ext^0(\ad E, \co(-p)) = H^0(\ad E(-p))$.  Since $E$
reduces to $T$, $\ad E$ splits as a sum of line bundles:
\begin{equation}
\label{s}
\ad E \cong \co^r \oplus \bigoplus_{\al \in \Phi} L_\al, 
\end{equation}
where $r$
is the rank of $G$ and $\Phi$ is the set of roots.  Of course
$H^0(\co^r(-p)) = 0$.  It is easily seen that $L_\al \cong \co(0 \,|\,
\al \pair \be_1 \,|\, \dots \,|\, \al \pair \be_n \,|\, 0)$ and hence
$L_\al(-p) \cong \co(-1 \,|\, \al \pair \be_1 \,|\, \dots \,|\, \al \pair
\be_n \,|\, 1)$.  Because of the weights $\pm 1$ at $p_\pm$, any
section of $L_\al(-p)$ vanishes on the first and last components.  By
\re{c}(c), $L_\al(-p)$ has a nonzero section if
and only if $\al \pair \be_i > 0 > \al \pair \be_j$ for some $i < j$.  On the
other hand, if $\al$ is a root, then so is $-\al$.  So $H^0(\bigoplus
L_\al(-p)) = 0$ if and only if, for all $\al$ and all $i,j$, $\al \pair
\be_i$ and $\al \pair \be_j$ are not of opposite sign.  This is equivalent
to all $\be_i$ lying in a common Weyl chamber.

(b): Consider the long exact sequence 
$$ 0 \lrow T^0_E \lrow T^0_{C_n,E} \lrow T^0_{C_n} \lrow T^1_E
\lrow T^1_{C_n,E} \lrow T^1_{C_n} \lrow 0.$$ The 1-dimensional
subspace $S \subset T^0_{C_n,E}$ arising from the $\G_m$-action
injects into $T^0_{C_n}$.  Hence $\dim T^0_{C_n, E} = 1$ if and
only if $T^0_E = 0$ and $T^0_{C_n}/S \to T^1_E$ is injective.  The
former is covered by (a), so it remains to consider the latter.

The connecting homomorphism in the sequence above may be described as
follows.  Let $D = \Spec k[\epsilon]/(\epsilon^2)$.  For $v \in
T^0_{C_n}$, an automorphism of $C_n \times_k D$ is given by $\id
+ \epsilon v$.  The pullback of $E \times_k D$ by this automorphism is
a bundle chain over $D$ whose isomorphism class determines an element of
$T^1_E$.

It is tempting to infer that the connecting homomorphism is zero.
After all, a first-order deformation cannot change the weights of the
$\G_m$-action on a line bundle, so line bundles and hence $T$-bundles
on $[C_n/\G_m]$ are rigid, but any $G$-bundle $E$ reduces to a
$T$-bundle, so the connecting homomorphism for $G$ factors through the
trivial one for $T$.  However, this ignores the framing.  The
deformation class of the unframed bundle is indeed zero in $H^1(\ad
E)$ by the above argument; but the framing at $p$ may deform
nontrivially.  From the long exact sequence of
$$ 0 \lrow \ad E(-p) \lrow \ad E \lrow \ad E \otimes \co_p \lrow 0,$$
namely \smallskip
\begin{equation}
\label{gg}
0 \lrow T^0_E \lrow H^0(\ad E) \lrow \g \oplus \g \lrow T^1_E \lrow
H^1(\ad E) \lrow 0,
\end{equation} 
it follows that the connecting homomorphism
lifts to $\g \oplus \g$.  It will be injective if and only if the
image of its lift intersects that of $H^0(\ad E)$ trivially.

A basis for $T^0_{C_n}$ consists of the infinitesimal generators $e_i$
for the $\G_m$-actions on each component of $C_n$.  The action on
the $i$th component lifts to act on the line bundle $\co(0 \,|\,
b_1 \,|\, \dots \,|\, b_n \,|\, 0)$ with weights $b_{i-1}$ and
$b_i$ on the two adjacent nodes.  Since the action on all other
components is trivial, this same lifting has weights $b_{i-1}$ and
$b_i$ on $p_+$ and $p_-$, respectively.

A $T$-bundle is essentially an $r$-tuple of line bundles, so the same
reasoning applies to the $T$-bundle to which $E$ reduces.  This has
weights $0, \be_1, \dots, \be_n, 0 \in \La$.  Hence the image of $e_i$
in $\g \oplus \g$ is $(\be_{i-1}, \be_i) \in \t \oplus \t \subset \g
\oplus \g$.

As for the image of $H^0(\ad E)$, it is now clear that only the part
meeting $\t \oplus \t$ is relevant.  Since this comes from $\co^r$ in
the splitting \re{s}, it is simply the diagonal
$\Delta_{\t} \subset \t \oplus \t$.  The composite map $T^0_{C_n}
\to (\t \oplus \t)/\Delta_{\t} \cong \t$ is then given by $e_i \mapsto
\be_i - \be_{i-1}$.

Since $S$ is spanned by $e_1 + \cdots + e_n$, this defines an
injection $T^0_{C_n}/S \to (\t \oplus \t)/\Delta_{\t}$ if and only
if $\be_1, \dots, \be_n$ are linearly independent in $\t = \La \otimes
k$.  

(c): There is a short exact sequence
$$1 \lrow \Auto E \lrow \Auto (C_n,E) \lrow \Auto C_n,$$
where as usual $\Auto E$ denotes the framed, $\G_m$-equivariant
automorphisms of $E$ lying over the identity on $C_n$.

The standard action of $\G_m$ on $C_n$ lifts to $E$, so there is
always a subgroup $\G_m \subset \Auto (C_n,E)$ lying over the
diagonal in $\Auto C_n \cong \G_m^{n+1}$.  On the other hand, the
Lie algebra of $\Auto E$ is $T^0_E$, and $\Auto E$ is connected 
\cite[6.7]{kansas}.  Hence $\Auto E$ is trivial if and only if $T^0_E = 0$.  Using
(a), if the $\be_i$ do not lie in a common Weyl chamber, then
certainly $\dim \Auto (C_n,E) > 1$, while if they do, then
$\Auto(C_n,E)$ is a subgroup of $\Auto C_n \cong \G_m^{n+1}$.
It suffices to show that, in the latter case, the $\be_i$ are
dependent in $\La$ if and only if $\Auto (C_n,E)/\G_m$ contains a
torus.  For any positive-dimensional subgroup of a torus contains a
torus \cite[IV 1.1.7]{demazure-gabriel}.

Indeed, we will show that the 1-parameter subgroup $\la: \G_m \to
\G_m^{n+1}$ given by $t^\la \deq (t^{a_1}, \dots, t^{a_n})$ for $a_i \in \Z$ lifts to
$\Auto (C_n,E)$ if and only if $\sum_{i=0}^n a_i(\be_{i+1}-\be_i) =
0$.  So there is always the diagonal $(1,\dots,1)$, but there will be
further subgroups if and only if the $\be_i$ are dependent.

Let $S \subset \Auto E$ be the 2-dimensional torus generated by $\la$
and the diagonal.  A lifting of $\la$ to $E$ is the same as a
$G$-bundle over $[C_n/S]$ whose pullback to $[C_n/\G_m]$ is $E$.  The
structure group of such a bundle reduces to $T$ \cite[6.4]{kansas}.  It is
elementary that a line bundle $\co(b_0 \st \cdots \st b_{n+1})$ admits
a $\la$-action preserving the framing at $p_\pm$ if and only if
$\sum_{i=0}^n a_i(b_{i+1}-b_i) = 0$.  Consequently, the $T$-bundle
$\F(\be)$ admits a $\la$-action preserving the framing if and only if
$\sum_{i=0}^n a_i(\be_{i+1}-\be_i) = 0$, as desired.  \fp

\bs{Corollary} 
\label{p}
Condition\/ {\rm (a)} is equivalent to $H^1(\ad E) = 0$ and to $\Auto E = 1$. 
\es  

\pf.  
Since the dualizing sheaf is $\om \cong \co(-p)$, by Serre duality
$H^1(\ad E)^* \cong T^0_E = 0$.  This is the Lie algebra of $\Auto E$.
On the other hand, $\Auto E$ is connected \cite[6.7]{kansas}, and a
connected group scheme over a field is trivial if and only if its Lie
algebra is trivial.  
\fp

\bs{Corollary} 
\label{hh}
If\/ condition {\rm (c)} holds, then $\Auto(C_n,E)/\G_m$ is a
sub-group scheme of $\G_m^n$, finite over $k$.
\es  

Hence in characteristic zero it is simply a finite group. 

\pf.  It is 0-dimensional and lies in the exact sequence
$$1 \lrow \Auto E \lrow \Auto(C_n,E)/\G_m \lrow \Auto C_n/\G_m.$$ But
$\Auto E = 1$ by \re{p}, so the last arrow {\red has trivial
scheme-theoretic kernel, hence} is a closed immersion
\cite[II 5.5.1b]{demazure-gabriel}, and $\Auto C_n / \G_m \cong
\G_m^n$.  
\fp

\bs{Corollary} 
For\/ {\rm (b)} to hold is an open condition in families, and likewise
for {\rm (c)}.  
\es

\pf.
By noetherian reduction, every bundle chain is obtained by base
change from one with a locally noetherian base.  The statement
for (b) then follows immediately from the usual semicontinuity
\cite[III 7.7.5]{ega}.  As for (c), the group scheme $\Auto(C,E)$ is of
finite type over its base, say by \re{ff} below.  Hence the dimension
at the identity section of $\Auto(C,E)$ is semicontinuous \cite[IV 13.1.3]{ega}.
\fp

\bs{Proposition} 
\label{t}
For $\be =(\be_1, \cdots, \be_n)\in \La^n$, let $L
\subset G$ be the Levi subgroup centralizing every $\be_i$, and let
$U_+, U_- \subset G$ be the unipotent subgroups consisting of root
spaces with some $\al \pair \be_i < 0$ (resp.\ $> 0$).  Then evaluation
at $p_\pm$ embeds the automorphism group $\Ga(\Ad E(\be))$ in $G
\times G$ with image $\Delta_L \ltimes (U_+ \times U_-)$.  
\es

\pf.  By definition, the kernel of the evaluation map is $\Auto
E(\be)$, which is trivial by \re{p}.  Hence the evaluation is a closed
immersion \cite[II 5.5.1b]{demazure-gabriel}.  It suffices, then, to
show that the image of the derivative is the Lie algebra of the
subgroup stated.  This follows from the splitting \re{s} together with
the observation that if $\al$ is a root of the Levi factor then the
line bundle $L_\al$ is trivial, whereas, thanks to \re{c}(c), if $\al$
is a root of $U_\pm$, then $L_\al$ has a section nonvanishing at
$p_\pm$ but vanishing at $p_\mp$. \fp

\bit{The stack of bundle chains}

Let $\BB$ denote the category of $G$-bundle chains {\red such that
$\G_m$ acts trivially over $p$ and} $H^1(\ad E) = 0$ on each geometric
fiber.  As we have seen, for a chain over $\Spec K$, these conditions
are equivalent to having {\red the isomorphism class of $E(\be)$,
where} $\be = (\be_1, \dots, \be_n) \in \La^n$ {\red and} $\be_i$ lie
in a common Weyl chamber.  For each such sequence $\be$ of weights,
the versal family $\E(\be) \to \C_n$ of \re{g} is an object of $\BB$
over $\A^n$.

\bs{Theorem}
\label{ff}
The category $\BB$ is a smooth Artin stack containing $BG \times
B\G_m$ as the dense open locus where {\red the length of the chain is}
$n=1$.  An atlas {\red is given by
$\bigsqcup_{\be} \E(\be): \bigsqcup_{\be} \A^{n(\be)} \to \BB$, where
$\be$ runs over all finite sequences in $\La$ lying in a common Weyl
chamber.}
\es

\pf.  We must prove three things: (1) that $\BB$ is a stack, namely
(a) a category fibered in groupoids, such that (b) descent data are
effective and (c) automorphisms are a sheaf; (2) that the diagonal
$\BB \to \BB \times \BB$ is (a) representable, (b) separated, and (c)
finitely presented; and (3) that the atlas stated is surjective and
smooth, so that the stack is algebraic.

Step (1).  That all morphisms over the identity are isomorphisms is
implicit in the definition of the category.  Also, pullbacks exist and
are unique up to isomorphism.  Hence $\BB$ is a category fibered in
groupoids, so (a) is true.  

Next, suppose an \'etale morphism $\tilde X \to X$ and descent data
for a bundle chain $\tilde E \to \tilde C \to \tilde X$ are given.
Since a bundle chain is equipped with actions of $\G_m$ on $\tilde E$
and $\tilde C$, it is understood that the descent data preserve these
actions.  Since algebraic spaces of finite type are a stack
\cite[04UD]{stacks}, $\tilde E \to \tilde C$ descend to algebraic
spaces $E \to C \to X$, and the $\G_m$-actions descend as well.  By
\re{h} $C$ is a chain over $X$.  A principal $G$-bundle such as
$\tilde E$ is by definition locally trivial in the \'etale topology,
so the same is automatically true of $E \to C$.  Hence (b) is true.

To see that automorphisms are a sheaf, observe that an automorphism of
a bundle chain $E \to C \to X$ is an automorphism of the total space
$E$ subject to the closed condition that it preserve the projection to
$X$ and commute with the action of $\G_m \times G$.  Again since
algebraic spaces of finite type are a stack, automorphisms of such
algebraic spaces are a sheaf, and hence the same is true of
automorphisms of a bundle chain.  This proves (c).

Step (2).  
Let $E_1, E_2$ be bundle chains over $Y$, corresponding to a morphism
$Y \to \BB \times \BB$.  It suffices to show that $\Hom(E_1,E_2)$ is
represented by a scheme $X$, separated and finitely presented over
$Y$.  The families $C_1, C_2$ of chains underlying $E_1, E_2$ are
induced by two morphisms $Y \to \CC$.  Hence by \re{h}
$\Hom(C_1,C_2)$ is represented by an algebraic space $Z$, separated
and finitely presented over $Y$, with a tautological morphism $g: C_1
\times_Y Z \to C_2$.  The $G$-bundle (or bitorsor) $\scrHom (E_1, g^*
E_2)$ over $C_1 \times_Y Z$ is locally trivial in the \'etale
topology, hence separated and finitely presented over $C_1
\times_Y Z$.  Isomorphisms $E_1 \to E_2$ are naturally equivalent
to sections of this bundle, and as such are represented by an
algebraic space, separated and locally of finite presentation over
$C_1 \times_Y Z$ \cite[6.2]{artin}, \cite[Appendix]{artin-versal}.

All that remains for Step (2) is to show that the diagonal is
quasi-compact.  This is more easily accomplished after Step (3).

Step (3).  Surjectivity follows from \re{p}, so it remains only to prove smoothness.

Let $E$ be a family of bundle chains over $Y$, corresponding to a
morphism $Y \to \BB$.  It suffices to show that $Y \times_\BB \A^n \to
Y$ is smooth.  By noetherian reduction we may assume that $Y$ is
locally noetherian.  Since smoothness descends under \'etale covers
\cite[1.15]{vistoli} we may also assume that $Y$ is a scheme.  By Step
(2) $Y \times_\BB \A^n$ is locally of finite presentation over $Y
\times \A^n$ and hence over $Y$.  So it suffices to show \cite[IV
17.5.4]{ega} that, for a local Artinian $k$-algebra $A$, any ideal $J
\subset A$, and any square
\begin{center}
\begin{tikzpicture}[>=angle 60] 
\matrix (m) [matrix of math nodes, row sep=2em, 
column sep=1em, text height=1.5ex, text depth=0.25ex] 
{ \Spec A/J & &  Y \times_\BB \A^n \\ 
\Spec A & & Y,\\ }; 
\path[->] 
(m-1-1) edge (m-1-3);
\path[->] 
(m-1-1) edge (m-2-1); 
\path[->] 
(m-1-3) edge (m-2-3); 
\path[->]
(m-2-1) edge (m-2-3);
\path[dashed,->]
(m-2-1) edge (m-1-3);
\end{tikzpicture}
\end{center}
there exists a diagonal arrow making the diagram commute.  For such an
$A$, we know $\m^n = 0$ for some $n > 0$, and $J/\m J,\m J/\m^2
J,\dots, \m^{n-1}J$ are finite-dimensional vector spaces over $A/\m$.
It therefore suffices to assume $\m J =0$.  The map $\Spec A \to Y$
gives a family of bundle chains over $A$.  

Another such family may be obtained as follows. Since $\A^n$ is smooth
and hence formally smooth over $k$, the morphism $\Spec A/J \to \A^n$
extends to $\Spec A \to \A^n$ by definition \cite[IV 17.1.1,
17.3.1]{ega}.  Pull back to $\Spec A$ the versal family $\E(\be)$ of
bundle chains over $\A^n$ from \re{g}.

We now have two families of bundle chains over $A$ with an isomorphism
of their restrictions to $A/J$.  It suffices to show that they are
isomorphic over $A$.  However, exactly as in Hartshorne
\cite[6.2a]{hartshorne2} or Sernesi \cite[1.2.12]{sernesi}, the set of
isomorphism classes of liftings of the bundle chain from $A/J$ to $A$
is a torsor for $H^1(\ad E) \otimes_{A/\m} J = 0$.  Hence the two
liftings are isomorphic.  This completes Step (3).

To finish Step (2), since $\bigsqcup_n \A^n \to \BB$ is surjective and
smooth, {\red base changing by
$\bigsqcup_m \A^m \times \bigsqcup_n \A^n \to \BB \times \BB$,} it
suffices to show that $\A^m \times_\BB \A^n \to \A^m \times
\A^n$ is quasi-compact \cite[1.15]{vistoli}.  Since $\A^m \times \A^n$
is affine, it suffices to show that $\A^m \times_\BB \A^n$ is
quasi-compact \cite[I 6.6.1]{ega}.  Express $\A^m$ as a disjoint union
over all $S \subset \{ 1, \dots m \}$ of $\A^m_S$, the locus where
$x_i = 0$ if and only if $i \in S$.  Then $\E(\be)|_{\A^m_S}$ is a
constant family, so
$\A^m_S \times_\BB \A^n_T \cong \A^m_S \times \A^n_T \times H$ where
$H$ (if nonempty) is the automorphism group of a bundle chain over
$\Spec k$.  The latter is an affine variety over
$k$ \cite[6.7]{kansas}.  Hence $\A^m_S \times_\BB \A^n_T$ is
quasi-compact.  But $\A^m
\times_\BB \A^n$ is the image of a surjective continuous morphism from
$\bigsqcup_{S,T} \A^m_S \times_\BB \A^n_T$,
so it is quasi-compact.
\fp

With the stacks $\CC$ and $\BB$ in hand, it is now easy to prove the
following.

\bs{Proposition}
\label{n}
Let $\E \deq \E(\be) \to \C_n$ as in {\rm \re{g}}.  Then the stack $H_\be$ of
automorphisms of $\E$ lying over the identity on $\C_n$ is
represented {\red over $\A^n$} by a 
smooth closed subgroup scheme of $\A^n \times G \times G$. 
\es

\pf.  To establish that $H_\be$ is represented by a smooth scheme,
it suffices to show that the natural morphism $\A^n \times_{\BB} \A^n
\to \A^n \times_{\CC} \A^n$ is smooth in a neighborhood of the
diagonal $\A^n \to \A^n \times_{\CC} \A^n$; for $H_\be \to \A^n$ is
nothing but the base change by this diagonal.

Since $\A^n \to \BB$ and $\A^n \to \CC$ are components of the atlases for
$\BB$ and $\CC$, they are smooth; hence $\A^n \times_\BB
\A^n$ and $\A^n \times_\CC \A^n$ are smooth over $\A^n$, and hence
over $k$.  Also, since the diagonals of $\BB$ and $\CC$ are separated
and finitely presented, $\A^n \times_\BB \A^n$ and $\A^n \times_\CC
\A^n$ are separated and finitely presented over $\A^n \times \A^n$ and
hence over $k$. In other words, they are smooth varieties.

The dimension of $\CC$ is $-1$, so the obvious morphism $\A^n \times
\G_m^{n+1} \to \A^n \times_\CC \A^n$ coming from the
$\G_m^{n+1}$-action on $\C_n$ is an open immersion.  The image of the
diagonal $\A^n \to \A^n \times_\CC \A^n$ lies in this open set.

Since the action of $\G_m^{n+1}$ on $\C_n$ lifts to $\E$, there are
compatible actions of $\G_m^{n+1}$ on $\A^n \times_{\BB} \A^n$ and
$\A^n \times_{\CC} \A^n$.  On the open set $\A^n \times \G_m^{n+1}$
mentioned above, the action is by multiplication on the second
factor.  

Since we have a morphism of smooth varieties, to show it is smooth
over $\A^n \times \G_m^{n+1}$, it suffices to show that its derivative
is everywhere surjective.  The compatible group actions show that the
image of the derivative contains the tangent space to $\G_m^{n+1}$.
But the smoothness of $\A^n \times_\BB \A^n$ over $\A^n$ shows that
the image of the derivative also contains a complement to this tangent
space.

Consequently, $H_\be$ is represented by a smooth variety, indeed a
group scheme smooth over $\A^n$.  Now let us show that this variety
admits a closed immersion (of group schemes) into $\A^n \times
G \times G$.  The immersion is given by restricting automorphisms of
$\E$ to the endpoints $p_\pm(\A^n)$, at which $\E$ is given the
conventional framings of \re{g}.  This defines the morphism
$H_\be \to \A^n \times G \times G$.  To show it is a closed immersion,
it suffices {\red \cite[II 5.5.1b]{demazure-gabriel}} to prove that it
{\red has trivial scheme-theoretic kernel}, {\red hence \cite[II
5.1.4]{demazure-gabriel} that the kernel is connected with trivial Lie
algebra.  That the Lie algebra is trivial follows from \re{gg}, while
connectedness is proved in another paper of the
authors \cite[6.7]{kansas}.}
\fp

\br{Remark}
\label{l}
Since $H_\be$ is a smooth family of groups over $\A^n$, with
connected fibers \cite{kansas}, it is connected.  Over the locus
$U \subset \A^n$ where $x_1 \cdots x_n \neq 0$, it is easy to
identify.  For $\E$ is trivial there, as explained in \re{g}.  Hence
$H_\be|_U \cong U \times G$.  Its immersion into $\A^n \times G
\times G$, determined by the conventional framings of \re{g}, is $(x, g) \mapsto
(x, g, x^{-\be} g\, x^\be)$, where $x^\be = \prod x_i^{\be_i}$ as in
\re{g}. 
\er

To get a Deligne-Mumford or tame Artin stack, we need to adjust $\BB$
in three ways: 
(1)~rigidify to get rid of the uniform $\G_m$-action;
(2)~pass to the $G\times G$-bundle parametrizing bundle chains equipped with
framings at $p_\pm$; and
(3)~restrict to the open substack where \re{j}(c) holds.
Of these, (3) is self-explanatory.  
Also, (2) is straightforward since $\BB$ clearly admits a universal
family $\pi: C \to \BB$ with sections $p_\pm:\BB \to C$ and a universal
$G$-bundle $E \to C$; the stack we want is nothing but $p_+^*E
\times_\BB p_-^*E$.  It remains to explain (1), which runs along the
lines laid down by Abramovich-Corti-Vistoli \cite{acv}.

Since the standard $\G_m$-action on chains acts on all of
our $G$-bundles by hypothesis, we may rigidify our moduli problem to
remove the common $\G_m$-symmetry.  This may be accomplished, as in the
theory of stable vector bundles, by introducing an equivalence on
families slightly broader than isomorphism.  Strictly speaking we
need to take the stack associated to this prestack \cite[\S5.1]{acv}.

Given a chain $C \to X$ and a line bundle $L \to X$, define a new
chain, denoted $C \otimes L$, to be the geometric quotient $(C
\times_B L^\times)/\G_m$.  That this satisfies the requirements of
being a chain over $X$ is easy to check: indeed, over each local
trivialization of $L$, it is isomorphic to $C$.  Thanks to the
$\G_m$-equivariance, there is a canonical equivalence between bundle
chains over $C$ and those over $C \otimes L$.  The desired prestack is
the one parametrizing bundle chains up to this equivalence.

Let $\MM$ be the stack obtained from $\BB$ by performing operations
(1), (2), (3) above.  It is the stack associated to the prestack
para\-met\-rizing bundle chains, framed at $p_\pm$, satisfying
\re{j}(c), modulo the equivalence described above.  It is easily
verified that the $\G_m$-action on $\MM$ induced by the standard
$\G_m$-action is trivial.  The same is true, by the way, of the stack $\NN$ of
rigidified chains obtained from $\BB$ by performing operation (1), and
there is a natural morphism $\MM \to \NN$.

By definition the {\it inertia stack} is $I \MM \deq \MM \times_{\MM
  \times \MM} \MM$.  It is naturally isomorphic to the stack whose
objects are pairs consisting of an object of $\MM$ and an automorphism
of that object, and whose morphisms are the obvious ones.  The two
projections induce naturally 2-isomorphic maps $I \MM \to \MM$, which
we regard as the same.

\bs{Theorem}
\label{k}
The stack $\MM$ has finite inertia: that is, $I\MM \to \MM$ is finite.
\es

\pf.  The fiber of $I\MM \to \MM$ over a bundle chain $E \to C \to X$
is represented by the scheme $\Auto (C,E)$.  We must show it is finite
over $X$.  That it is quasi-finite follows immediately from \re{j}(c).
It therefore suffices to show that $I\MM \to \MM$ is proper.  This we
accomplish with the valuative criterion.

Since $\MM$ has a locally noetherian atlas (whose connected components
are $\A^n \times G \times G$) it is locally noetherian by definition
\cite[4.7.1]{lmb}.  Since it is algebraic, the diagonal $\MM \to \MM
\times \MM$ is representable and finitely presented, hence so is
its base change $I \MM \to \MM$.  According to
Laumon-Moret-Bailly \cite[7.8, 7.12]{lmb}, it suffices to show that
for any complete discrete valuation ring $R$ over $k$ with fraction
field $F$, any commutative diagram in the form of the right-hand square of   
\begin{center}
\begin{tikzpicture}[>=angle 60] 
\matrix (m) [matrix of math nodes, row sep=2em, 
column sep=1em, text height=1.5ex, text depth=0.25ex] 
{ \Spec \tilde{F} & & \Spec F & & I\MM \\ 
\Spec \tilde{R} & & \Spec R & & \MM\\ }; 
\path[->] 
(m-1-1) edge (m-1-3);
\path[->] 
(m-1-3) edge (m-1-5);
\path[->]
(m-2-1) edge (m-2-3);
\path[->] 
(m-2-3) edge (m-2-5);
\path[->] 
(m-1-1) edge (m-2-1); 
\path[->] 
(m-1-3) edge (m-2-3); 
\path[->] 
(m-1-5) edge (m-2-5); 
\path[dashed,->]
(m-2-1) edge (m-1-5);
\end{tikzpicture}
\end{center}
may be extended by adding a left-hand square (where $\tilde{R}$ is a
discrete valuation ring containing $R$ whose residue field $\tilde{F}$
is an extension of $F$) and a morphism as indicated by the dashed
arrow shown, unique up to 2-isomorphism.  By the Cohen structure
theorem \cite[I 5.5]{hartshorne} $R \cong K[[t]]$ and $F
\cong K((t))$ for some extension $K/k$.  Since $K[[t]] = \lim_{n \to
  \infty} K[t]/(t^n)$, the morphism $\Spec R \to \MM$ may be lifted to
$f: \Spec R \to \A^n \times G \times G$, a connected component of the
smooth atlas for $\MM$.  Let $\E(\be) \to \C_n \to \A^n$ be the
corresponding bundle chain.  Then $E \deq f^* \E(\be)$ is a framed bundle
chain over $K[[t]]$; we must show that any given automorphism of
$E_\eta \deq E|_{K((t))}$ extends uniquely to all of $E$.  

Let $i_1, \dots, i_m$ denote those coordinates of $f_1: \Spec R \to
\A^n$ which vanish identically.  There may be no such coordinates, but
if there are, then the chain $C \deq f^* \C_n$ will be reducible.  In
any case, its generic fiber $C_\eta$ consists of $m+1$ lines, and $C$
itself consists of $C = C^{(1)} \cup \cdots \cup C^{(m+1)}$ glued end
to end.  

The given automorphism of $E_\eta$ is a $K((t))$-valued point of
$\Auto(C_n,E)/\G_m$, which by \re{hh} is diagonalizable and
0-dimensional, so it is of finite order, say $n$, prime to the
characteristic.  Consequently, $E_\eta$ together with its given
automorphism may be regarded as a bundle chain over
$[C_{m+1}/(\G_m \times \mu_n)]$.

Over some extension $\tilde{K}$ of $K$, the structure group of such an
object reduces to $T$.  The proof of this is the same as \cite[4.3,
6.4]{kansas}, with two exceptions.  Part B of 4.3 is replaced by an
easy argument: having a $\G_m$-action trivializes any $G$-bundle over
the dense $\G_m$-orbit in $[\P^1/\G_m]$, and the $\mu_n$-action {\red
there, being essentially a homomorphism $\mu_n \to G$, may then} be
conjugated into $T$ \cite[7.1]{kansas}, so the associated $G/B$-bundle
has a section over this dense orbit, which extends to a regular
section over $[\P^1/\G_m]$ by the valuative criterion.  And in 6.4,
the centralizer $Z(\chi(\G_m \times \mu_n))$ need not be connected,
but {\red as the centralizer of a torus in a reductive group,} its
identity component is reductive \cite[A.8.12]{cgp}, so, over some
extension $\tilde{K}$, it has a Bruhat decomposition of the form $BWB$
just like a (connected) reductive group.

Consequently, once $E_\eta$ is reduced to $T$ as a bundle chain over
$[C_{m+1}/\G_m]$, the given automorphism is uniquely determined by its
triviality over $p_\pm$ together with the automorphism of
$C_\eta \cong C_{m+1}$ underlying it.

Suppose the latter is given by the $m+1$-tuple $(z_0, \dots, z_m)$ for
$z_i \in K((t))^\times$.  Each $z_i$ extends to a uniform action on
$C^{(i)}$ in the sense of
\re{i}.  Indeed, each $C^{(i)}$ is the base change of a versal chain
by $\pi f_1: \Spec R \to \A^{k_j}$, where $k_j = i_{j+1} - i_j -1$ and
$\pi: \A^n \to \A^{k_j}$ is projection on the coordinates between
$i_j$ and $i_{j+1}$.  Let $\F^{(i)}$ be the $T$-bundle over $C_{k_j}$
associated to $\L_0^\times \times \cdots \times \L_{k_j+1}^\times$ via
the homomorphism $(\be_{i_j}, \dots, \be_{i_{j+1}}): \G_m^{k_j+2} \to
T$, and let $\E^{(i)}$ be the associated $G$-bundle.  Then $(\pi
f_1)^* \E^{(i)} \cong E|_{C^{(i)}}$.  Since the uniform action of
$\G_m$ on $\A^{k_j}$ lifts to $\F^{(i)}$ and thence to $\E^{(i)}$, the
uniform action of $z_i$ on $C^{(i)}$ lifts to $E|_{C^{(i)}}$.

These actions have compatible weights at the endpoints of $C^{(i)}$.
So they glue to give an automorphism of $E$, which reduces to $T$
and hence extends the one given over $E_\eta$.  Two such
automorphisms agree on $E_\eta$ and hence everywhere, giving the
required uniqueness.  \fp

An Artin stack with finite inertia is defined by
Abramovich-Olsson-Vistoli \cite{aov} to be {\it tame\/} if the
push-forward of quasi-coherent sheaves to the coarse moduli space is
an exact functor.  They prove \cite[3.2]{aov} that an Artin stack with
finite inertia is tame if and only if every geometric point has
linearly reductive stabilizer.  In particular, an Artin stack over
a field of characteristic zero is tame if and only if it is
Deligne-Mumford.

\bs{Theorem} The stack $\MM$ is tame, and the open substack
where condition\/ {\rm \re{j}(b)} holds is Deligne-Mumford.
\es

\pf.  If $E \to C$ is a bundle chain over a field, then 
there is a
short exact sequence of group schemes
$$1 \lrow \Auto E \lrow \Auto (C,E) \lrow \Auto C.$$ But $\Auto E$ is
trivial by \re{p}, so $\Auto (C,E) \to \Auto C$ is a closed
immersion \cite[II 5.5.1b]{demazure-gabriel}.  Hence $\Auto(C,E)$ is a
sub-group scheme of a torus, so it is diagonalizable and hence
linearly reductive.

For a stack with finite inertia to be Deligne-Mumford, it suffices
\cite{aov} that every geometric point have reduced stabilizer.  But
that is exactly what is guaranteed by \re{j}(b).  \fp

\bs{Proposition} 
Neither $\BB$ nor $\MM$ is separated unless $G=1$.
\es

\pf.  First let $G = \G_m$, so that $\La = \Z$.  Let $x$ be the
coordinate on $\A^1$, and consider the chain $\C_1 = \Bl(\A^1 \times
\P^1, 0 \times [0,1])$ over $\A^1$.  For any $i\in \Z$, the functions
$1$ and $x^i$ provide framings, at $p_+$ and $p_-$ respectively, for
the $\G_m$-bundle $\L_1^i = \O(-iE)$, where $E$ is the exceptional
divisor of the blow-up.  Let $\L$ be the bundle chain $\L_1^2$ with
this framing, but let $\L'$ be the base change of $\L_1$ by $x \mapsto
x^2$ with the induced framing.  Then the nonzero fibers of $\L$ and
$\L'$ are isomorphic as framed bundle chains away from 0, but the
fibers at 0 are not: rather, they are $\O(0|2|0)$ and $\O(0|1|0)$
respectively.

Consequently, $(\L, \L')$ defines a morphism $\A^1 \to \MM \times \MM$
such that $\A^1 \times_{\MM \times \MM} \MM = \A^1 \sans 0$, which is
not proper over $\A^1$.  For $G = \G_m$, consequently, $\MM$ is not
proper over $\MM \times \MM$, and hence $\MM$ is not separated.
Exactly the same families define morphisms to $\BB$ and show that it
is not separated either.

For general reductive $G \neq 1$, let $\be: \G_m \to G$ be a nontrivial
cocharacter.  Extension of structure group by $\be$ converts
the $\G_m$-example above into a $G$-example showing that $\BB$ and $\MM$
are not separated in this case either.  \fp

{\red

Changing the framings at $p_\pm$ defines an action of $G \times G$ on
$\MM$.  Each orbit is the locus where the isomorphism class of the
unframed bundle chain is fixed.  For example, the orbit of $E
= \O_{\P^1}$ has isotropy group $\{(g,g^{-1}) \st g \in G\}$ and may
therefore be identified with $G$.  Clearly this action lifts to the
obvious action on $\A^{n(\be)} \times G \times G$, the component of
the atlas for $\MM$ obtained from the versal chain $\E(\be)$.

\bs{Proposition}
Every $G \times G$-orbit closure in $\MM$ is smooth, and every
intersection of two orbit closures is in a normal crossing.  In
particular, the complement $\MM \sans G$ of the dense orbit $G$ is a
divisor with normal crossings.
\es

\pf.
The same is true of the preimages of the orbits in each component
$\A^{n(\be)} \times G \times G$ of the atlas for $\MM$, since these
preimages are nothing but the orbits of the usual action of
$\G_m^{n(\be)}$ on $\A^{n(\be)}$.
\fp
}

\bit{The stability condition}

We now cut down our stack to a separated stack of finite presentation
by imposing some sort of stability condition.  This requires that we make
some non-canonical choices.

First, as in toric geometry, choose a fan $\Si$, supported within
$\La_\Q^+$, consisting of rational simplicial cones.  Also choose a
nonzero lattice element $\be_i \in \La$ on each ray of $\Si$.  {\red
The choice of a distinguished lattice point on each ray is precisely the additional
structure making a fan into a {\em torsion-free stacky fan} as defined
by Borisov-Chen-Smith \cite{bcs}, who established a correspondence
between such objects and toric orbifolds.  (The slightly more general
notion introduced by Borisov-Chen-Smith of a {\em stacky fan},
possibly with torsion, will not be needed here.)}

Second, choose an ordering $\be_1,\dots, \be_N \in \La$ of these
elements, and hence of all the rays of $\Si$.  This induces an
ordering of the rays of any $\si \in \Si$.  Thus to any
$n$-dimensional $\si \in \Si$ is associated an increasing function, by
abuse of notation also denoted $\si: \{1,\dots,n\} \to \{1,\dots,N\}$,
such that $\si = \langle \be_{\si(1)}, \dots, \be_{\si(n)} \rangle$.

\br{Definition}
\label{m}
A bundle chain over a field is said to be $\Si$-{\em
  stable} if it is isomorphic to the bundle chain $E(\be_{\si(1)}, \dots,
\be_{\si(n)})$ for some $n$-dimensional cone $\si \in \Si$.  
An arbitrary bundle chain is $\Si$-{\em stable}
if its geometric fibers are $\Si$-stable. 
\er

Strictly speaking, this depends not only on the stacky fan $\Si$, but
also on the ordering of its rays.  Notice that a $\Si$-stable bundle
chain automatically satisfies \re{j}(c), hence belongs to $\MM$.

\bs{Theorem} 
\label{r}
The locus of $\Si$-stable bundle chains is an open {\red $G \times
G$-invariant} substack ${\red \MM_G}(\Si) \subset \MM$, of finite presentation
over $k$, {\red which} contains $G$ as a dense open substack.  It is always
separated and is proper if and only if the support of $\Si$ equals
$\La_\Q^+$.
\es

\pf.  For any $\Si$-stable $E(\be) \to C_n$ over $k$, the versal chain
$\E(\be)$ over $\A^n \times G \times G$ is also $\Si$-stable.  The
$\Si$-stable locus is therefore open.

Finitely many such families now constitute an atlas, one for each cone
in $\Si$.  The stack is therefore of finite presentation over $k$.

It is clear from this atlas that the smaller locus where the chain has
no nodes is a dense open substack.  There the bundle must be trivial,
as it is acted on with weight 0 over $p_\pm$.  A global trivialization
is determined by the framing at $p_+$, relative to which the framing
at $p_-$ induces a natural isomorphism of this substack with $G$.


{\red It remains to establish separation and properness.  This will be
accomplished by reducing to a split maximal torus $T \subset G$ and
invoking the corresponding facts about toric stacks.  

To prove separation, it suffices \cite[7.8]{lmb} to show, for}
any complete discrete valuation ring
$R$ over $k$ with fraction field $F$ and any two morphisms $f_1, f_2:
\Spec R \to {\red \MM_G}(\Si)$, that any 2-isomorphism $f_1|_{\Spec F} \to
f_2|_{\Spec F}$ extends to a 2-isomorphism $f_1 \to f_2$.  That is, if
$E_1$ and $E_2$ are the framed bundle chains over $\Spec R$ induced by
$f_1$ and $f_2$, {\red then any isomorphism $E_1 \cong E_2$ over
$\Spec F$ extends over $\Spec R$.}

As in the proof of \re{k}, $R \cong K[[t]]$ and $F \cong K((t))$ for
some extension field $K/k$.  Moreover, following Deligne-Mumford
\cite[4.19]{dm}, \cite[II 7.2.3]{ega}, one may assume that the
image $f_1(\eta) \cong f_2(\eta)$ of $\eta = \Spec K((t))$ lies in the
dense open substack $G \subset {\red \MM_G}(\Si)$.

{\red For any} $f:
\Spec K((t)) \to G$, there exists an unique 1-parameter subgroup $\la:
\Spec K((t)) \to T$ in the positive Weyl chamber, and morphisms
$\ga,\de: \Spec K[[t]] \to G$, such that $f = \ga \la \de^{-1}$.
{\red This was originally proved by Iwahori-Matsumoto \cite{im}.  It
is best regarded as being essentially}
the Bruhat decomposition for the formal
loop group $G((t))$ whose points are morphisms $\Spec K((t)) \to G$.
{\red To see this, let} ${\mathcal G} = G((t)) \rtimes \G_m$, let ${\mathcal B}$ be the
inverse image of the Borel $B$ under the evaluation map $G[[t]]
\rtimes \G_m \to G$, let ${\mathcal N} = \La \rtimes (N(T) \times
\G_m)$, and let ${\mathcal S}$ be the usual generators for the affine
Weyl group $\tilde W = \La \rtimes W$.  Then $\mathcal{G, B, N, S}$
satisfy \cite[6.2.8]{kumar} the requirements to be a Tits system
\cite[5.1.1]{kumar}, and, since $\La^+ = W \backslash \tilde W / W$,
we obtain \cite[5.1.3c]{kumar} the {\red Bruhat decomposition}
$$G((t)) = \bigsqcup_{\la \in \La^+} G[[t]] \,\, \la \,\, G[[t]].$$

Taking $f = f_1$ in the above, and then replacing $f_1$ and $f_2$ by
$\ga^{-1}f_1\de$ and $\ga^{-1}f_2\de$, we may assume that $f_i|_{\Spec
  K((t))}$ is a 1-parameter subgroup $\la$ of $T$, indeed in the
positive Weyl chamber.  

{\red We claim that in fact $f_1$ and $f_2$ factor through
$\MM_T(\Si)$, the moduli stack with structure group $T$.  To establish
this, we will first lift them to the atlas for $\MM_G(\Si)$, then
adjust the lifts to go into the atlas for $\MM_T(\Si)$.}

The connected components of the atlas for ${\red \MM_G}(\Si)$ are of
the form $\A^n \times G \times G$ {\red for $n = n(\be)$}.  Choose a
point $(x(0),g(0),h(0))$ in one of them lying over $f_1(0) \in
{\red \MM_G}(\Si)$.  By a suitable choice of component we may assume
that $x(0) = 0$.

Since $K[[t]]$ is a limit of Artinian $K$-algebras, we may lift $f_1$
to $\tilde f_1 = (x,g,h): \Spec K[[t]] \to \A^n \times G \times G$.
Then $\tilde f_1(\eta) \subset U \times G \times G$, where $U \subset
\A^n$ is the locus where $x_1 \cdots x_n \neq 0$, for this is the
inverse image of $G \subset {\red \MM_G}(\Si)$.  Indeed, using the conventional
framings of \re{g} we see that the morphism $U \times G \times G \to
G$ is $(x,g,h) \mapsto g^{-1} x^\be h$. 

As pointed out in \re{g}, the action of $\G_m^n$ on $\E(\be) \to \C_n
\to \A^n$ fixes the conventional framing at $p_+$ but changes it at
$p_-$.  Hence there is an action of $\G_m^n$ on $\A^n \times G \times
G$ that projects to the obvious action on the first factor, acts
trivially on the second factor, acts nontrivially on the third factor,
and lifts to the framed bundle chain $\E(\be)$.  {\red It therefore
preserves the morphism to $\MM_G(\Si)$.}

{\red Hence} we may adjust the lift $\tilde f_1$ by acting by any
morphism $z: \Spec K[[t]] \to \G_m^n$.  That is to say, $z
\tilde f_1$ is another lift of $f_1$.  We may {\red thus} assume that the
first factor of $\tilde f_1$, namely $x: \Spec K((t)) \to U$, is a
1-parameter subgroup $t \mapsto (t^{i_1},
\dots, t^{i_n})$, with each $i_j > 0$, when $U$ is
identified with $\G_m^n$.  

We then have $\la = g^{-1}\, x^\be\, h$ where both $\la$
and $x^\be$ are 1-parameter subgroups in the positive Weyl chamber
$\La^+$.  By the uniqueness in the Iwahori decomposition, it follows
that $x^\be =\la$.  Hence $h = x^\be\, g\, x^{-\be}$.

Comparison with \re{l} now reveals that $(g,h)$ is a regular section
over $\Spec K[[t]]$ of $\tilde f_1^*H_\be \subset \Spec K[[t]] \times
G \times G$. Acting by the inverse of this automorphism takes the
framing to a trivial framing.  Hence the pullbacks of $\E(\be)$ by the
morphisms $\tilde f_1 = (x,g,h)$ and $(x,1,1)$ are isomorphic as
framed bundle chains over $\Spec K[[t]]$.  The latter is therefore
also a lift of $f_1$.

In other words, we may assume that $\tilde f_1 = (x,1,1)$ for some
1-parameter subgroup $x:\Spec K[[t]] \to \G_m^n$.  Likewise, we may
assume a similar lift for $f_2$.  Thus the structure group of $E_1$
and $E_2$, even with their framings, is reduced to $T$.  {\red This
establishes our claim that} $f_1$ and $f_2$ factor through
$\MM_T(\Si)$.

But the latter is easily seen to be the toric stack $\X(\Si)$ {\red
associated to the stacky fan $\Si$ by} Borisov-Chen-Smith \cite{bcs}.
Indeed, the atlas consists of components $\A^n \times T \times T$, but
by \re{l}, the automorphism group $H_\be$ reduces, when $G=T$, to the
constant family consisting of the diagonal $\Delta_T \subset T \times
T$.  Since this acts by framed bundle automorphisms, we may descend to
the quotient $(\A^n \times T \times T)/ \Delta_T
{\red \, \cong \,} \A^n \times T$.  This is essentially the toric atlas of
Borisov-Chen-Smith \cite[4.4]{bcs}: indeed, there is a further
automorphism group $\G_m^n$, and $[(\A^n \times T)/\G_m^n] = [\A^n /
K]$, where $K = \ker \be:\G_m^n \to T$ is {\red a torus times} the
finite group called $N(\si)$ by Borisov-Chen-Smith.  (In positive
characteristic this is a zero-dimensional sub-group scheme of
$\G_m^n$.)

Since $\X(\Si) = \MM_T(\Si)$ is tame, the natural morphism to
$X(\Si)$, its coarse moduli space \cite[3.7]{bcs}, is
proper \cite[1.1]{conrad}.  Since $X(\Si)$, {\red being a toric
variety,} is certainly separated \cite[I \S2]{kkms}, it finally
follows that $f_1
\cong f_2$ in $\MM_T(\Si)$ and hence in ${\red \MM_G}(\Si)$.  Hence the
latter is separated.  

To prove properness, proceed as in the proof of separation and
observe that the limits of all 1-parameter subgroups in the positive Weyl
chamber $\La^+$ exist in $X(\Si)$, and hence in $\X(\Si)$ which is
proper over it, if and only if the support of $\Si$ is $\La_\Q^+$.
\fp

A connected atlas for ${\red \MM_G}(\Si)$ may be exhibited as follows.  Let
$\be \deq (\be_1, \dots, \be_N)$, and let $\E(\be)$ be the
corresponding versal bundle chain.  Also let $\A_\be \deq \A^N$, and
let $\A_\be \times G \times G$ parametrize all framings at $p_\pm$ of
bundle chains in $\E(\be)$ by mapping $(x,g_+,g_-) \mapsto \Psi_\pm
\circ g_\pm$, where $\Psi_\pm: \A_\be \times G \to p_\pm^* \E(\be)$
are the conventional framings of \re{g}.  The pullback $\tilde\E(\be)$
of $\E(\be)$ to $\A_\be \times G \times G$ is then tautologically
framed.

While $\tilde\E(\be)$ is not $\Si$-stable everywhere, it is
$\Si$-stable over an open subset.  As in \re{d}(c), for any $S \subset
\{ 1, \dots, N\}$, let $U_S = \{ x \in \A^N \st x_i \neq 0 \mbox{ if }
i \not \in S\}$.  Let $\A^0_\Si$ be the union of all $U_S$ such that
$\langle \be_i \st i \in S \rangle$ is a cone in $\Si$.  This
precisely covers the permitted splitting types, so the restriction of
$\tilde\E(\be)$ to $\A^0_\Si \times G \times G$ is $\Si$-stable.

\bs{Proposition}
\label{o}
The morphism $\A^0_\Si \times G \times G \to {\red \MM_G}(\Si)$ defined by the
framed bundle chain $\tilde\E(\be)$ is smooth and surjective,
providing a connected atlas for ${\red \MM_G}(\Si)$.
\es

\pf.  Surjectivity is clear, since the bundle chains that appear are
precisely those permitted by $\Si$-stability.  As for smoothness,
observe that by \re{d}(c), the projection $\pi_S: U_S \to \A^{|S|}$
induces an isomorphism $\C_N |_{U_S} \cong \pi_S^* \C_{|S|}$; since
the line bundles $\L_i$ correspond under this isomorphism, it follows
that $\E(\be)|_{U_S} \cong \pi_S^* \E(\be_S)$, where $\be_S$ is the
subsequence of $\be$ indexed by $S$.  The restriction of $\A^0_\Si
\times G \times G \to {\red \MM_G}(\Si)$ to $U_S \times G \times G$ therefore
factors as a composition of $\pi_S$ with the morphism $\A^{|S|} \times
G \times G \to {\red \MM_G}(\Si)$ provided by the obvious atlas for ${\red \MM_G}(\Si)$,
both of which are smooth.  \fp

If the structure group $G$ is a torus, then ${\red \MM_G}(\Si)$ is a toric
stack. Indeed, we have seen in the proof of \re{r} that it is the toric
stack $\X(\Si)$ associated to the stacky fan $\Si$ by
Borisov-Chen-Smith \cite{bcs}, {\red generalizing} Cox \cite{cox}.  Moreover,
the connected atlas of \re{o} coincides in that case with their global
quotient construction.  

\bs{Theorem}
\label{cc}
The closure of $T$ in $\MM_G(\Si)$ is $\overline T = \X(W\Si)$, the toric stack
associated to the Weyl-invariant stacky fan $W\Si$.
\es

\pf. First, $\MM_T(W\Si) = \X(W\Si)$ by the remark above.

Second, note that since $E(w\be) \cong E(\be)$ as $G$-bundles, the
framed $T$-bundle chains parametrized by $\MM_T(W\Si)$ are all
$\Si$-stable $G$-bundle chains, determining a representable morphism
$\MM_T(W\Si) \to \MM_G(\Si)$.

Third, since a proper, unramified, and universally injective morphism
is a closed embedding \cite[02K5]{stacks}, it suffices to show that
this representable morphism enjoys these three properties.

Properness is immediate when the support of $\Si$ is all of
$\La_\Q^+$, for then both $\MM_T(W\Si)$ and $\MM_G(\Si)$ are proper
over $k$ by \re{r}.  For general $\Si$, extend $\Si$ to a stacky
simplicial fan $\overline{\Si}$ whose support is $\La_\Q^+$ and
observe that by \re{b}, the inverse image in $\MM_T(W\overline{\Si})$ of
$\MM_G(\Si) \subset \M_G(\overline{\Si})$ is $\MM_T(W\Si)$.

Since $\MM_T(W\Si)$ and $\MM_G(\Si)$ are {\red separated and} finitely
presented over $k$, our morphism is finitely presented \cite[IV
1.6.2]{ega}, hence is unramified if and only if it is formally
unramified.  It therefore suffices, as in the proof of \re{h}, to show
that the derivative is everywhere injective.  For this we recall the
description of the deformation space in terms of the sheaf ${\mathcal
S}$ in \re{q}.  Let $F \to C$ be a framed $T$-bundle chain over a
geometric point of $\MM_T(W\Si)$, and let $E \to C$ be the associated
$G$-bundle chain in $\MM_G(\Si)$.  Let ${\mathcal S}_G$ and ${\mathcal
S}_T$ denote the sheaves on $C$ of $G$-invariant (resp.\
$T$-invariant) K\"ahler differentials on $E$ (resp.\ $F$).  Then there
is a short exact sequence
$$0 \lrow \bigoplus_{\al \in \Phi} L_\al \lrow {\mathcal S}_G \lrow {\mathcal
  S}_T \lrow 0,$$ where $\Phi$ denotes the roots of $G$ and $L_\al$ is
the line bundle defined in the proof of \re{q}.  By that result, the relevant
derivative is the induced map $\Ext^1_C({\mathcal S}_T,\O(-p))
\to \Ext^1_C({\mathcal S}_G,\O(-p))$.  Since $H^0(L_\al(-p)) = 0$ as
shown in the proof of \re{q}, it follows that the derivative is injective.

Finally, to see that the morphism is universally injective, let $F,F'
\to C_n$ be two $W\Si$-stable framed $T$-bundles over a field, and
suppose their extensions $E,E'$ to $G$ are isomorphic as framed
$G$-bundles.  Let $g: E \to E'$ be the isomorphism taking framings to
framings.  By \re{b}, the splitting types of $E$ and $E'$ agree modulo
the Weyl group.  That is, there is an isomorphism $f$ of the
associated $N(T)$-bundles, though not necessarily compatible with the
framings.  The automorphism $f^{-1} \circ g$ of $E$ then preserves, at
the points $p_\pm$, the $N(T)$-structure induced by $F$.  Hence it
must lie in $\Delta_{L \cap N(T)} \subset \Delta_L$ by \re{t}, where
$L$ is the Levi subgroup described there.  Since $L$ centralizes each
$\be_i$, it is simply the automorphism globally induced by a fixed
element of $L \cap N(T)$.  The structure group of $f^{-1} \circ g$,
and hence of $g$, thus reduces to $N(T)$ everywhere on the chain
$C_n$, meaning that $g$, viewed as an isomorphism from the total space
of $E$ to that of $E'$, takes the total space of the $N(T)$-bundle
containing $F$ to the one containing $F'$.  The total spaces of these
$N(T)$-bundles contain $F$ and $F'$ as connected components, and we
know that $g(F) = F'$ since the framings lie there and are preserved
by $g$.  Hence the structure group of $g$ reduces to $T$, that is, $F
\cong F'$ as framed bundles.  \fp

\bit{The Vinberg monoid}

We now turn to a different description of our moduli stacks: as global
quotients, by a torus action, of an open subset in an algebraic
monoid.  As an application, we deduce in \S7 necessary and
sufficient conditions for their coarse moduli spaces to be projective.

Our description generalizes Vinberg's description \cite{v} of the
wonderful compactification of an adjoint group as a quotient, by a
torus action, of an open subset in a certain algebraic monoid.
Vinberg calls it the {\it enveloping semigroup}, and we call it the
{\it Vinberg monoid}.  The construction of this monoid has been
extended to reductive groups by Alexeev and Brion \cite{ab}.  At least
some parts have been extended to positive characteristic by Rittatore
\cite{rit1, rit2}, using the theory of spherical varieties.

At this point we must assume that our ground field $k$ is algebraically
closed.  All of the above references assume this, and the most
comprehensive account, namely Vinberg's, also assumes that the
characteristic is zero.  We think it plausible that Vinberg's results
listed below could be established for a split reductive group over an
arbitrary field, using the work of Huruguen \cite{huruguen} on
spherical varieties over arbitrary fields.  If so, our subsequent
results would go through.

To state Vinberg's results, we first define some polyhedral cones.
Let $\V$ be the character lattice of $T \subset G$.  Let the Weyl
group $W$ act on $\V_\Q \oplus \V_\Q$ by $w(\mu,\nu) \deq
(\mu,w\nu)$.  Let the simple roots $\al_i$ and fundamental weights
$\varpi_j$ be indexed by $i,j \in \Om$, the Dynkin diagram.  Define
$K$ to be the polyhedral cone in $\V_\Q \oplus \V_\Q$ spanned by the
simple roots in $\V_\Q \oplus 0$ and the diagonal, that is, $K =
\langle (\al_i,0) \st i \in \Om \rangle + \Delta_{\V_\Q}$.  Its
intersection $K^+ \deq K \cap (\V_\Q^{\phantom{+}} \oplus \V_\Q^+)$
with the positive Weyl chamber is the product of a simplicial cone,
spanned by $(0, -\al_i)$ and $(\varpi_j,\varpi_j)$, with a linear
subspace.  The union of the images of
$K^+$ under the action of $W$ is again a cone: indeed $WK^+ =
\bigcap_{w \in W} wK$ \cite[5.11]{renner}.

Vinberg relates the faces of $WK^+$ to those of $K^+$ as follows.  For
$I,J \subset \Om$, define cones
$$D_I = \langle \al_i \st i \in I \rangle,
\qquad 
C_J = \langle \varpi_j \st j \in J \rangle,$$ 
\begin{equation}
\label{fij}
F_{I,J} = \{ (\mu, \nu) \in \V_\Q \oplus \V_\Q
\st \mu -\nu \in D_I,\, \nu \in C_J + \Delta_{\V_\Q}^W \}.
\end{equation} 
The $F_{I,J}$ are the faces of $K^+$.  A face $F_{I,J}$ is said to be
{\it essential\/} if no connected component of $\Om \sans J$ is
contained in $I$ when these are regarded as subsets of the Dynkin
diagram.  Vinberg proves that every face $F$ of $WK^+$ satisfies $WF =
WF_{I,J}$ for exactly one essential face $F_{I,J}$ of $K^+$.

\br{Definition}
\label{jj}
Let $Z$ be a torus with a given isomorphism to the maximal torus $T$
of $G$.  The {\it Vinberg monoid} is an affine algebraic monoid $S_G$ over
$k$ with group of units $(Z \times G)/Z_G$, where the center $Z_G$ of
$G$ is antidiagonally included in $Z \times G$.  {\red It is defined
to be the spectrum of the ring of matrix entries of
representations of $(Z \times G)/Z_G$ whose highest weights lie
in $WK^+$.}  
\er
 
{\red The Vinberg monoid} is equipped with
monoid homomorphism{\red s} $\pi: S_G \to \A$ and $\psi:\A \to S_G$
{\red satisfying $\pi \circ \psi = \id_\A$}.  Here $\A \deq \Spec
k[\al_i]$ is an affine space having $Z/Z_G$ as its group of units.
{\red Vinberg \cite{v} proves the following:}

\smallskip

\noindent \vv{1} 
The restriction of $\pi$ to the group of units is the projection to
the torus $Z/Z_G$, and the restriction of $\psi$ to $Z/Z_G$ is given by
$\psi(z) = (z,z^{-1})$.

\noindent \vv{2} 
The closure in $S_G$ of $(Z \times T)/Z_G$ is the affine toric
variety $X((WK^+)^\vee)$. 

\noindent \vv{3}
Each orbit of $Z \times G \times G$ in $S_G$ (acting by left and
right multiplication) contains exactly one orbit of $Z \times N(T)$ in
$X((WK^+)^\vee)$.  The orbits $O_{I,J} \subset S_G$ are therefore
indexed by the essential faces.  

\noindent \vv{4}
The morphism $U_- \times Z \times \A \times U_+ \to S_G$ given by
$(g_-,z,x,g_+) \mapsto g_- \, z \, \psi(x) g_+$, where $U_\pm$ are the usual
maximal unipotents in $G$, is an open immersion with image $S^0_G \deq
\bigcup_I O_{I,\Om}$.  

\noindent \vv{5} 
This open subset has as geometric quotient $S^0_G/Z \cong G_{\ad}$,
the wonderful compactification of the adjoint group $G_{\ad} \deq
G/Z_G$.

\smallskip

In fact, this geometric quotient is actually a geometric invariant
theory quotient, as we shall now prove.  Since $S_G$ is affine, we may
choose as our ample line bundle the trivial bundle $\O$, but with the
usual action twisted by a character $\rho \in \V$ of $Z$.  The
geometric invariant theory quotient is then $\Proj \bigoplus_{i \geq
  0} k[S_G]^{i\rho}$, where $k[S_G]^{i\rho}$ is the space of regular
functions with weight $i \rho$ for the usual action.  In this setting,
the Hilbert-Mumford numerical criterion may be stated as follows.

\bs{Lemma}
\label{w}
If a torus acts effectively on an affine variety, linearized on
a trivial bundle by a character $\rho$, then $x$ is semistable (resp.\
stable) if and only if for every nontrivial cocharacter $\la$ such
that $\lim_{t \to 0} t^\la \, x$ exists, $\rho \cdot \la \geq 0$ (resp.\
$> 0$).
\es

A proof of this version of the criterion is given by King
\cite{king} and by Gulbrandsen-Halle-Hulek \cite{ghh}.

\bs{Theorem}\label{x}
When linearized by a character $\rho$ in the interior of the positive
Weyl chamber $\V^+$, the $Z$-action on $\SG$ has semistable (and
stable) set $\SG^0$, and the geometric invariant theory quotient $\SG
\, \modmod_{\hspace{-0.2em}\rho} \, Z$ is the wonderful
compactification $\overline{G}_{\ad}$.  \es

\pf. Since the stable and semistable sets are unions of orbits of $Z
\times G \times G$, by \vv{2} and \vv{3} it suffices to
consider points in the affine toric variety $X((WK^+)^\vee)$, the
closure of the maximal torus $(Z\times T)/ Z_G$.  In particular, to
every essential face $F_{I,J}$ of $K^+$ corresponds a Weyl orbit of
faces of $WK^+$.

We will test for stability for the $Z$-action using the numerical
criterion \re{w}. 
If $x$ lies in the torus orbit corresponding to a face $F$ of $WK^+$,
it is straightforward to check that $\lim_{t \to 0} t^\la \, x$
exists if and only if $F \cdot \la \geq 0$, meaning $\mu \cdot
\la \geq 0$ for all $\mu \in F$. Hence the torus orbit is stable
(resp.\ semistable) if and only if for all nontrivial $\la \in \La$ with $F \cdot
\la \geq 0$, one also has $\rho \cdot \la > 0$ (resp.\ $\geq 0$).

Now $Z$, its subgroups $\la$, and its character $\rho$ are acted on
trivially by $W$, whereas $WF = WF_{I,J}$ for every face $F$ of
$WK^+$.  So the stability condition above holds for $F$ if and only if
it holds for $F_{I,J}$.  Hence $O_{I,J}$ is semistable (resp.\ stable)
if and only if for all nontrivial $\la \in \La$ with $F_{I,J} \cdot
\la \geq 0$, one also has $\rho \cdot \la \geq 0$ (resp.\ $> 0$).

Observe from \re{fij} that $\mu \in D_I + C_J$ if and only if $(\mu,
\nu) \in F_{I,J}$ for some $\nu$.  For any 1-parameter
subgroup $\la$ of the first factor $Z$ of $Z \times T$, then, $F_{I,J}
\cdot \la = (D_I + C_J) \cdot \la$.  

If $J = \Omega$, it follows that $F_{I,J} \cdot \la \geq 0$ implies
$\rho \cdot \la > 0$, for $\rho$ is in the interior of the positive
Weyl chamber $\V^+$.  Consequently, each $O_{I,\Om}$, and hence all of
$\SG^0$, is stable.

If $J \neq \Omega$, however, then from the definition of essential
pair there is necessarily some $i \in \Omega$ contained in neither $I$
nor $J$.  In this case $-\al_i^\vee$ provides a destabilizing
1-parameter subgroup, since $\varpi_i \cdot \al_i^\vee = 1$ whereas
for $j \neq i$, $\alpha_j \cdot \alpha_i^{\vee}\leq 0$ and $\varpi_j
\cdot \alpha_i^{\vee}=0$.  Hence the complement of $\SG^0$ is
unstable.

The geometric invariant theory quotient $\SG \,
\modmod_{\hspace{-0.2em}\rho} \, Z$ is therefore a geometric quotient
$\SG^0/Z$, and as such agrees with $\overline{G}_{\ad}$ by
\vv{5}.

As an alternative to the last step, we may assume that $\rho$ is a
regular dominant weight of $G_{\ad}$, and consider the representation
$\overline{R}^\rho:\SG \to \End V_\rho$ also defined by Vinberg.  His
work implies \cite[(57),(59)]{v} that $\overline{R}^\rho(O_{I,J})=0$
if and only if $J\neq \Omega$.  It immediately follows that
$\SG^0=(\overline{R}^\rho)^{-1}(\End V_\rho\setminus \{0\})$.  We
therefore get an embedding of $G/Z_G=G_{\ad}$ into
$\mathbb{P}\left(\End V_\rho \right)$, which may be regarded as the
geometric invariant theory quotient $\End V_\rho \,
\modmod_{\hspace{-0.2em}\rho} \, Z$.  Hence $\SG^0/Z$ will contain the
closure of $G_{\ad}$ in $\mathbb{P}\left(\End V_\rho \right)$, which
is well known to be the wonderful compactification of $G_{\ad}$
\cite[2.1]{ej}.  Since both have exactly $2^r$ orbits of
$G_{\ad}\times G_{\ad}$, where $r$ is the rank of $G_{\ad}$, they are
isomorphic. \fp

\bit{The Cox-Vinberg hybrid}

We aim to hybridize the Vinberg quotient, which realizes
$\overline{G}_{\ad}$ as a torus quotient of the monoid $\SG$, with the
Cox construction, which realizes toric orbifolds as torus quotients of
affine spaces $\A_\be$ \cite{bcs, fmn}.  The hybrid will realize
${\red \MM_G}(\Si)$ as a torus quotient of a monoid $S_{G,\be}$ obtained as the
base change of $\SG$ by a linear map of affine spaces $\A_\be \to \A$.
This was inspired by Brion's observation \cite{brion1} that spectra of
Cox rings of spherical varieties, such as the wonderful
compactification, are often base changes of the Vinberg monoid.

We briefly recall the Cox construction.  Let $Z$ be a torus and $\La$
its cocharacter lattice.  In $\La_\Q$, let $\Si$ be a torsion-free
stacky {\red simplicial} fan, that is, a {\red simplicial} fan
equipped with a set of nonzero lattice elements $\be_1, \dots, \be_N$
so that each ray of the fan contains exactly one $\be_i$.  Let
$\A_{\be}\deq \A^N$ and $\G_\be \deq
\G_m^N$.  For any $\si \subset \{ 1, \dots, N\}$, let $U_\si \deq \{x \in
\A_\be \st x_i \neq 0 \mbox{ if }i \not\in \si \}$.  Then let $\A_\be^0$
be the union of all $U_\si$ such that  $\langle \be_i \st i \in
\si\rangle$ is a cone in $\Si$.  Though it is not explicitly indicated in
the notation, this open subset $\A_\be^0 \subset \A_\be^{\phantom{0}}$ depends on
$\Si$ as well as $\be$.

Let $\phi_\be:\G_\be \to Z$ be given by $\phi_\be(z_1,\dots, z_N) \deq
\prod z_i^{\be_i}$ (or more concisely, $\phi_\be(z) = z^\be$) and let
$K_\be \deq \ker \phi_\be$.  If the $\be_i$ span $\La_\Q$, then
\begin{equation}\label{short2}1\to K_{\be} \to \G_\be
  \stackrel{\phi_\be}{\longrightarrow} Z\to 1.\end{equation} 
In this case Borisov-Chen-Smith \cite{bcs} show that 
$\X(\Si) \deq [\A_\be^0/ K_{\be}^{\phantom{0}}]$ is a separated tame 
stack whose coarse moduli space is the toric variety $X(\Si)$.
(They consider only characteristic zero, but this construction is general.) 

\br{Remark}
Obviously $[\A_\be^0 /
K_{\be}^{\phantom{0}}] \cong \left[ \left(\A_\be^0\times Z\right)
/ \G_\be^{\phantom{0}}\right]$ if the $\be_i$ span $\La_\Q$.  But even
if they do not, so that \re{short2} is not exact on the right,
$\left[ \left(\A_\be^0\times Z\right) / \G_\be^{\phantom{0}}\right]$
is still a toric stack $\X(\Si)$ agreeing with the stack of that name
constructed by Fantechi-Mann-Nironi \cite[7.12]{fmn}.
\er

Now let $Z$ be isomorphic to the maximal torus $T$ of $G$, and suppose
the support of $\Si$ lies in the positive Weyl chamber $\La_\Q^+$.
Let $\pi:Z \to Z/Z_G$ be the projection, and let $\overline{\phi}_\be
\deq \pi \circ \phi_\be: \G_\be \to Z/Z_G$.  Since $\be_i \in \La^+$,
the group homomorphisms $\overline{\phi}_{\be_i}: \G_m \to Z/Z_G$
extend to monoid homomorphisms $\A^1 \to \A$.  Hence
$\overline{\phi}_\be$ extends to a monoid homomorphism $\A_\be \to
\A$.

\br{Definition}
Let the {\it Cox-Vinberg monoid\/} be the fibered product $S_{G,\be} \deq \A_\be
\times_{\A} \SG$.  Likewise, let $S_{G,\be}^0 \deq \A_\be^0
\times_{\A} \SG^0$.
\er

The Cox-Vinberg monoid is a reductive monoid, flat over $\A_\be$, with
group of units $\G_\be\times G$.  Like the Vinberg monoid, it has a
projection $\pi_\be: S_{G,\be} \to \A_\be$ whose restriction to the
group of units is $(z,g) \mapsto z$ and a section $\psi_\be: \A_\be
\to S_{G,\be}$ whose restriction to $\G_\be \subset \A_\be$ is $z
\mapsto (z,z^{-\be})$.  By the way, choosing $\A_\be \cong \A$ is
generally not allowed: for then $\be_i$ must be (some permutation of)
the fundamental coweights, and these are not in $\La$ unless the
semisimple part of $G$ has trivial center.

\bs{Theorem}  
The stack $[S_{G,\be}^0/ \G_\be^{\phantom{0}}]$ is
canonically isomorphic to ${\red \MM_G}(\Si)$.
\label{u}
\es 

\pf.  Identify $\A_\be \times G \times G$ with the space parametrizing
all framings at $p_\pm$ of bundle chains in $\E(\be)$ via $(x,g_+,g_-)
\mapsto \Psi_\pm \circ g_\pm$, where $\Psi_\pm: \A_\be \times G \to
p_\pm^* \E(\be)$ are the conventional framings of \re{g}.  Then, over
the base $\A_\be$, the natural action of the group scheme $H_\be$
on $\A_\be \times G \times G$ is 
$${\red h} (x,g_+,g_-) \deq (x,e_+({\red h})g_+,e_-({\red h}) g_-),$$ 
where $e_\pm:H_\be \to G$ are the
evaluations at $p_\pm$ relative to $\Psi_\pm$.  The pullback
$\tilde\E(\be)$ of $\E(\be)$ to $\A_\be \times G \times G$ is then
tautologically framed, and the $H_\be$-action lifts to the
framed bundle chain $\tilde\E(\be)$.

As seen in \re{n}, $\pi \times e_+ \times e_-: H_\be \to \A_\be
\times G \times G$ is a closed immersion of smooth group schemes over
$\A_\be$.  Hence the quotient $J_\be \deq (\A_\be \times G \times G) /
H_\be$ is a smooth family of homogeneous spaces over $\A_\be$,
and in particular a smooth scheme.  Since the $H_\be$-action lifts to
$\tilde\E(\be)$, this framed bundle chain descends to the quotient.
This determines a morphism $J_\be \to \BB$.

Moreover, let $J_\be^0 \deq (\A^0_\be \times G \times G) /
H_\be$ be the open subset of $J_\be$
lying over the aforementioned $\A^0_\be \subset
\A_\be$.  There, the framed bundle chain is $\Si$-stable and hence
defines a morphism $J_\be^0 \to {\red \MM_G}(\Si)$.  This is surjective since
the atlas $\A^0_\be \times G \times G \to {\red \MM_G}(\Si)$ of \re{o} factors
through it.  This also implies that it is smooth.  For a morphism $g$
is smooth provided that $g$ is finitely presented, $f$ is smooth and
surjective, and $g \circ f$ is smooth \cite[02K5]{stacks}.  Apply this
to $\A^0_\be \times G \times G \to J_\be^0 \to {\red \MM_G}(\Si)$, observing
that $J_\be^0 \to {\red \MM_G}(\Si)$ is finitely presented since $J_\be^0 \to
\Spec k$ and the diagonal ${\red \MM_G}(\Si) \to {\red \MM_G}(\Si) \times {\red \MM_G}(\Si)$
are.  It follows that $J_\be^0 \to {\red \MM_G}(\Si)$ is an atlas.

The torus $\G_\be$ acts on every ingredient in this recipe, hence on
$J_\be^0$ and the framed bundle chain over it.  The morphism $J_\be^0
\to {\red \MM_G}(\Si)$ therefore descends to $[J_\be^0/\G_\be^{\phantom{0}}]
\to {\red \MM_G}(\Si)$.  We claim this is an isomorphism.  It suffices
\cite[3.8]{lmb} to show that the action morphism $J_\be^0 \times
\G_\be \to J_\be^0 \times_{{\red \MM_G}(\Si)} J_\be^0$ is an isomorphism.  
{\red For this it} suffices to show \cite[IV 17.9.1]{ega} that it
is \'etale and bijective on points over any field.  That it is
bijective is straightforward, since bundle chains of distinct
splitting types lie over the distinct $\G_\be$-orbits in $\A_\be$, and
the isomorphism classes of possible framings are bijectively
parametrized by the points of $J_\be^0$.

To show that it is \'etale, since source and target are smooth and
finitely presented, it suffices to show that its derivative is an
isomorphism on Zariski tangent spaces.  Let $E \to C$ be a bundle
chain in $J_\be^0$.  From the long exact sequence
of 
$$0 \lrow \ad E(-p) \lrow \ad E \lrow \g \oplus \g \lrow 0$$
and \re{p},
the tangent space to the
fiber of $J_\be^0 \to \A_\be$ is $H^1(\ad E(-p))$.
So there are short exact sequences
$$\begin{array}{ccccccccc}
0 & \lrow & 
H^1(\ad E(-p)) & \lrow & 
T_{C,E} \, J_\be^0 \oplus T_z \G_\be & \lrow & T_C \A_\be \oplus T_z \G_\be 
& \lrow & 0\\
& & \down & & \down & & \down && \\
0 & \lrow & \Delta_{H^1(\ad E(-p))} & \lrow &
T_{C,E}(J_\be^0 \times_{{\red \MM_G}(\Si)} J_\be^0 )
& \lrow &
T_{C\times z\cdot C} (\A_\be \! \times_\NN \! \A_\be) & \lrow & 0.
\end{array}
$$
The right-hand column is the derivative of the action $\A_\be
\times \G_\be \to \A_\be \times_\NN \A_\be$, where $\NN$ is the stack of rigidified
chains as in \S3.  It is easy to check, using $T^1_C = \Ext^1(\Omega,
\O(-p))$ as in \re{i} and \re{q}, that this derivative is an
isomorphism.  Hence the middle column is an isomorphism as well.  This
completes the proof that $[J_\be^0/\G_\be] \to {\red \MM_G}(\Si)$ is an
isomorphism.  

It suffices, finally, to exhibit a $\G_\be$-equivariant isomorphism
$J_\be^0 \to S^0_{G,\be}$.  Let $\G_\be$ act on $\A^0_\be \times G
\times G$ by 
$$z \cdot (x,g_+,g_-) = (zx,\, g_+,\, g_-z^{-\be}).$$  
Then the morphism $\A^0_\be \times G \times G \to S^0_{G,\be}$ is
given by $(x,g_+,g_-) \mapsto g_+ \psi_\be(x) g_-^{-1}$ is
$\G_\be$-equivariant.  Indeed, as a closed condition the equivariance
may be verified on the group of units where $\psi_\be(x) =
(x,x^{-\be})$.  This morphism is smooth with image $S^0_{G,\be}$ by
\vv{4}.  The inverse image of $\psi_\be^{\phantom{0}}(\A^0_\be)$ is
therefore a smooth family of subgroups of $G \times G$.  By \vv{1},
over a point $z \in \G_\be^{\phantom{0}} \subset \A^0_\be$, this group
is exactly $(1 \times z^\be)\, \Delta_G\, (1 \times z^{-\be})$.  By
the last remark in \re{g}, it therefore coincides with $H_\be$.

Consequently, this morphism descends to a $\G_\be$-equivariant,
surjective morphism $J_\be^0 \to S^0_{G,\be}$ of schemes over
$\A^0_\be$.  It is universally injective, since on each fiber of
$J_\be^0$ over $\A^0_\be$ we have exactly divided by the stabilizer
group of the transitive $G \times G$-action.  And it is \'etale, since
source and target are smooth of the same dimension, and the original
$\A^0_\be \times G \times G \to S^0_{G,\be}$ is smooth, so the
derivative is everywhere surjective.  Again by the fundamental
property of \'etale morphisms \cite[IV 17.9.1]{ega}, $J_\be^0
\to S^0_{G,\be}$ is an isomorphism.  \fp

\bs{Corollary} 
If $G$ has trivial center, $\Si$ comprises the single cone $\La_\Q^+$,
and $\be_i = \varpi_i^\vee$, the fundamental coweights, then the
coarse moduli space of ${\red \MM_G}(\La_\Q^+)$ is $M(\La_\Q^+) \cong
\overline{G}_{\ad}$, the wonderful compactification of
$\overline{G}_{\ad} \deq G/Z_G$.  \es

\pf.  In this case $\phi_\be$ is an isomorphism and so $S_{G,\be}
\cong S_G$ as varieties with the action of a torus $\G_\be \cong Z$.
By \vv{5}, $S_G^0/Z \cong \overline{G}_{\ad}$ as a geometric
quotient.  \fp

\br{Remark} 
Since $\A^0_\be$ and hence $S^0_{G,\be}$ do not depend on the choice
of ordering of $\be_1, \dots, \be_N$, by \re{u} the moduli stack
${\red \MM_G}(\Si)$ does not depend on it either.  In other words, for two
different orderings, there is a functor taking $\Si$-stable bundle
chains for one ordering to $\Si$-stable bundle chains for the other.
We do not know an explicit description of this functor.  \er

{\red

\bit{Functoriality}

We now turn to the question of functoriality: given a homomorphism
$f:G \to G'$ of split reductive groups, for which fans does it extend
to a morphism $\MM_G(\Si) \to \MM_{G'}(\Si')$ of the
compactifications?  The answer turns out to be essentially the same as
for toric stacks.


Let $\Si$, $\Si'$ be torsion-free stacky fans supported on the
positive Weyl chambers, with distinguished elements $\be_i$ and
$\be'_j$, respectively.  By composing with inner automorphisms,
we may suppose without loss of generality that $f(T) \subset T'$ and that $\langle f^*
\al'_j \rangle \subset \langle \al_i \rangle$.  Thus $f$ determines a
linear $f_*: \t \to \t'$ (well-defined modulo a subgroup of $W \times W'$).
Say that $f_*$ defines a {\em morphism of stacky fans} $W\Si \to
W'\Si'$ if for every cone $\si \in W\Si$ there exists $\si' \in
W'\Si'$ such that $f_*(\si) \subset \si'$ and if for every
distinguished $\be_i \in \si$, $f_*(\be_i)$ is an integer
combination of those $\be'_j \in \si'$.  (Because the fans are
simplicial, these integers are unique and nonnegative.)
Borisov-Chen-Smith show \cite[4.5]{bcs} that a homomorphism $T
\to T'$ of tori extends to a morphism $\X(\Si) \to \X(\Si')$ of
toric orbifolds if and only if its derivative $\t \to \t'$
defines a morphism of stacky fans $\Si \to \Si'$.

\bs{Theorem}
A homomorphism $f:G \to G'$ extends to a morphism of stacks
$\MM_G(\Si) \to \MM_{G'}(\Si')$ if and only if $f_*: \t \to \t'$
defines a morphism of stacky fans $W\Si \to W'\Si'$.
\es

\pf.
By \re{cc}, the closures of $T$ and $T'$ in $\MM_G(\Si)$ and
$\MM_{G'}(\Si')$ are $\X(W \Si)$ and $\X(W' \Si')$ respectively, so the
necessity is an immediate consequence of Borisov-Chen-Smith's
criterion \cite[4.5]{bcs}.

Conversely, suppose that $f_*$ defines a morphism of stacky fans.
Then each $f_* \be_i = \sum_j a_{ij} \be'_j$, where the integers
$a_{ij}$ satisfy $a_{ij} = 0$ unless $\be'_j \in \si' \supset
f_* \si \ni \be_i$.  Consider the monoid homomorphism
$\rho: \A_\be \to \A_{\be'}$ given by $\rho(x_i) = \left( \prod_i
x_i^{a_{ij}} \right)$.  Thanks to the vanishing condition just
mentioned, if $x_i \neq 0$ for $i \not\in \si$, then $\prod_i
x_i^{a_{ij}} \neq 0$ for $j \not\in \si'$.  That is,
$\rho(U_\si) \subset U_{\si'}$, where $U_\si$ and $U_{\si'}$ are the
open sets in the Cox construction.  Consequently
$\rho(\A^0_\be) \subset \A^0_{\be'}$, where as before $\A^0_\be
= \bigcup_{\si \in \Si} U_\si$.

Furthermore, since $\langle f^* \al'_j \rangle \subset \langle \al_i \rangle$, the
homomorphism $\bar{f}: Z/Z_G \to Z'/Z_{G'}$ induced by $f$ extends to
a monoid homomorphism $\bar{f}: \A \to \A'$.  The square
\begin{center}
\begin{tikzpicture}
\matrix (n)[matrix of math nodes,  row sep=3em, column sep=2.5em, text height=1.5ex, text depth=0.25ex]
{\A_\be & \A_{\be'} \\ \A & \A' \\};
\path[->] 
(n-1-1) edge[shorten >=2pt,shorten <=2pt] node[above] {$\rho$} (n-1-2);
\path[->]
(n-1-1)  edge[shorten >=2pt,shorten <=2pt] node[auto, swap] {$\overline{\phi}_\be$} (n-2-1);
\path[->] 
(n-2-1) edge[shorten >=2pt,shorten <=2pt] node[above] {$\bar{f}$} (n-2-2);
\path[->] 
(n-1-2) edge[shorten >=2pt,shorten <=2pt] node[auto] {$\overline{\phi}_{\be'}$} (n-2-2);
\end{tikzpicture}
\end{center}
%
commutes, since 
\beqas
\overline{\phi}_{\be'}(\rho(x_i)) 
& = & \overline{\phi}_{\be'}\left( \prod_i x_i^{a_{ij}} \right) \\
& = & \prod_{i,j} x_i^{a_{ij}\be'_j} \\
& = & \prod_i x_i^{\sum_j a_{ij}\be'_j} \\
& = & \prod_i x_i^{f_* \be_i} \\
& = & f\left( \prod_i x_i^{\be_i} \right) \\
& = & \bar{f}(\overline{\phi}_\be(x_i)).
\eeqas

On the Cox side, therefore, we have a morphism
$\A^0_\be \to \A^0_{\be'}$, equivariant for the actions of the groups
of units $\G_\be$ and $\G_{\be'}$, and lying over the monoid
homomorphism $\bar{f}$.

On the Vinberg side, it follows directly from the definition of
the Vinberg monoid \re{jj} that $f:G \to G'$ extends to a
homomorphism $S_f: S_G \to S_{G'}$.  By \vv{1} $S_f$ takes
$\psi(\A)$ to $\psi'(\A')$ and hence $S^0_G = G \psi(\A) G$ to
$S^0_{G'} = G' \psi'(\A') G'$.

There is hence a morphism of stacks $\rho \times S_f:
[(\A^0_\be \times_\A S^0_G)/\G_\be] \to [(\A^0_{\be'} \times_\A'
S^0_{G'})/\G_{\be'}]$, which by \re{u} is exactly
$\MM_G(\Si) \to \MM_{G'}(\Si')$.  Restricting to the open subset
$G \cong \G_\be \times_\A (G \times Z)/Z_G \subset \A^0_\be \times_\A
S^0_G$, we recover the original homomorphism $f: G \to G'$.  The
homomorphism therefore extends as desired.
\fp

Since the stacks in question represent moduli problems, there
must be a functor associated to the morphism of stacky fans
above, taking stable $G$-bundle chains to stable $G'$-bundle
chains.  Alas, we do not know an explicit description of this
functor.

}

\bit{\boldmath $M_G(\Si)$ as a geometric invariant
  theory quotient}

In this section we show that if the coarse moduli space of the toric
orbifold $\MM_T(W\Si)$ is projective, then the same is true of
$\MM_G(\Si)$, and the quotient construction of \re{u} is a geometric
invariant theory quotient.  In fact we show a little more.

Slightly adapting the terminology of Cox et al.\ \cite{coxbigbook} and
Hausel-Sturmfels \cite{haussturm}, we say a variety is {\it
  semiprojective} if it is projective over an affine variety.  In
particular, a proper semiprojective variety is projective.  Also, a
closed subvariety of a semiprojective variety is semiprojective.  A
variety is semiprojective if and only if it is Proj of an integral
algebra of finite type over $k$, so any geometric invariant theory
quotient of a semiprojective variety is semiprojective.

We say a rational polyhedral fan is a {\it normal fan} if it
comprises the normal cones of some polyhedron, or,
equivalently \cite[7.2.4, 7.2.9]{coxbigbook}, if its support is convex
and admits a strictly convex piecewise linear function to $\Q$, linear
on each cone of the fan.  Then a fan $\Si$ is a normal fan if and
only if the toric variety $X(\Si)$ is semiprojective.

\bs{Theorem}\label{v}
If $\Si$ is a normal fan and $W\Si$ has convex support, then the
coarse moduli space $M_G(\Si)$ of $\MM_G(\Si)$ is semiprojective and indeed
is a geometric invariant theory quotient $S_{G,\be}\modmod \G_\be$
with semistable (and stable) set $S^0_{G,\be}$.  
\es

\bs{Corollary}
If $\Si$ is a normal fan and has support equal to $\La_\Q^+$, then the
coarse moduli space $M_G(\Si)$ of $\MM_G(\Si)$ is projective.
\es

\pf.  It is proper by \re{r} and semiprojective by \re{v}.  \fp

Before proving the theorem, we pause to observe a purely polyhedral consequence.

\bs{Corollary}
If $\Si$ is a normal fan and its support is the entire Weyl chamber
$\La_\Q^+$, then $W\Si$ is a normal fan.
\es

\pf.  The support of $W\Si$ is then the entire $\La_\Q$, hence convex,
so $M_G(\Si)$ and hence $M_T(W\Si) \subset M_G(\Si)$ are
semiprojective varieties.  \fp

\bs{Lemma}\label{z}
If $W\Si$ has convex support, then the projection of any element in
$|\Si|$ onto any face of the positive Weyl chamber is also contained in
$|\Si|$.  \es
 
\pf. The projection of any $x \in |\Si|$ onto any face of the Weyl
chamber is contained in the convex hull of the Weyl orbit $Wx$, hence
belongs to the convex $W$-invariant set $|W\Si|$.  \fp

\pf\ of {\rm \re{v}}.  Since ${\red \MM_G}(\Si) \cong [S^0_{G,\be}/ G_\be]$ as
stacks, it suffices to find a linearization of the $\G_\be$-action on
$S_{G,\be}$ for which the semistable (and stable) set is
$S^0_{G,\be}$.  Since $S_{G,\be}$ is affine, we may choose the trivial
bundle as our ample bundle.  However, $\G_\be$ should act by a
nontrivial character.

We will combine two characters adapted to the two factors of
$S_{G,\be} = \A_\be \times_\A S_G$.  First, let $\xi: \G_\be \to \G_m$
be any character so that $(\A_\be \times
Z)\modmod_{\hspace{-0.2em}\xi} \, \G_\be = X(\Si)$ with semistable
(and stable) set $\A^0_\be \times Z$.  Note that this is equivalent to
$\A^0_\be$ being the stable set for the $K_\be$-action.  Such a
character $\xi$ always exists when $\Si$ is a normal fan
\cite[Remark 12.1, p.\ 202]{dolg}.  Second, let $\rho:\G_\be \to \G_m$
be the composition of $\phi_\be: \G_\be \to Z$ with a character $Z \to
\G_m$ in the interior of the positive Weyl chamber.  Then linearize
the $\G_m$-action on $S_{G,\be}$ ``asymptotically'' by $\xi + m \rho$
for $0 \ll m \in \Z$.  (Intuitively, the motivation is that $m \rho$
translates the moment polyhedron for the $K_\be$-action so that all of
its vertices lie in the positive Weyl chamber.)

By the numerical criterion \re{w}, $(x,y)$ is semistable (resp.\
stable) for the action of $\G_\be$ if and only if for every nontrivial
cocharacter $\la$ such that $\lim_{t \to 0} t^\la \, (x,y)$
exists, $(\xi + m \rho) \cdot \la \geq 0$ (resp.\ $>0$).

But suppose first that $(x,y) \in S^0_{G,\be} = \A^0_\be \times_\A
S^0_G$, which means that $(x,1) \in \A_\be \times Z$ and $y \in S_G$
are stable points for the $\G_\be$ actions linearized by $\xi$ and
$\rho$ respectively.  And suppose that $\lim_{t \to 0}t^\la \, (x,y)$
exists.  Then certainly $\lim_{t \to 0}t^\la \, y$ exists, so $\rho
\cdot \la > 0$.  However, $\lim_{t \to 0}t^\la \, (x,1)$ exists if and only if
$\phi_\be (t^\la) = 1$.  In this case $\rho \cdot \la = 0$ and $\la$
acts trivially on $S_G$, but since $(x,1)$ is stable, $\xi \cdot \la >
0$.  Hence $(\xi + m \rho) \cdot \la > 0$, as required.  Otherwise, if
$\phi_\be (t^\la) \neq 1$, then since $1 \ll m$ and $\rho \cdot \la
> 0$, again $(\xi + m \rho) \cdot \la > 0$.  Hence $(x,y)$ is a stable
point for the $Z$-action on $S_{G,\be}$ linearized by $\xi + m \rho$.

Now, conversely, suppose that $(x,y) \in S_{G,\be} = \A_\be \times_\A
S_G$ is stable for the linearization $\xi + m \rho$, and let $\la$ be
a nontrivial cocharacter of $\G_\be$.  If $\phi_\be (t^\la) = 1$, so
that the 1-parameter subgroup corresponding to $\la$ is contained in
$K_\be$, then $\rho \cdot \la = 0$ and $\la$ acts trivially on $S_G$,
so the linearization reduces to $\xi$ on $\A_\be$, showing that $x \in
\A^0_\be$, the stable set for the $K_\be$-action.

What remains is to show that in this case $y$ is also stable.  In fact
we will show, equivalently, that if $y$ is unstable then $(x,y)$ is
unstable.  The unstable sets in $S_G$ and $S_{G,\be}$ are unions of
orbits of $Z \times G \times G$ and $\G_\be \times G \times G$,
respectively.  By \vv{2} and \vv{3} it suffices to assume that $y$
lies in the affine toric variety $X((WK^+)^\vee) \subset S_G$, and
indeed in an orbit $O_{I,J}$ corresponding to an essential face of $WK^+$.  On
the other hand, the $\G_\be$-orbits of $\A_\be$ are
indexed by $H \subset \{ 1, \dots, N\}$ and are, of course, simply
$$\tilde{O}_H \deq \{ x \in \A_\be \st x_j \neq 0 \Leftrightarrow j \in
H\}.$$
Since the morphism $\A_\be \to \A$ is induced by contraction with $\be:\langle
\al_i \st i \in \Om \rangle \to \Q_{\geq 0}^N$, it follows that $\tilde{O}_H
\subset \A_\be$ maps to $O_I \subset \A$ where $$I=\left\{i\in\Omega \st
\alpha_i\!\cdot\!\beta_j =0 \mbox{ for all } j\not\in H\right\}.$$

Recall now from the proof of \re{x} that if $y \in O_{I,J}$ is
unstable for an essential pair $(I,J)$, then $J \neq \Om$, there
exists $i\not\in I \cup J$ by definition of an essential pair, and
$-\al_i^{\vee}: \G_m \to Z$ then proves to be a destabilizing
cocharacter for $O_{I,J}$.  To destabilize $(x,y) \in
\tilde{O}_H \times_{O_I} O_{I,J}$, by the numerical criterion \re{w}
it suffices to lift a multiple of $-\al_i^\vee$ to a 1-parameter
subgroup $\la: \G_m \to \G_\be$ such that $\lim_{t \to 0}t^\la \, x$
exists.  For the necessary inequality $(\xi + m \rho) \cdot \la < 0$
is then automatic from $1 \ll m$.

If $t^\la = (t^{\ell_1},\dots,t^{\ell_N})$ and $x \in \tilde{O}_H$, then
clearly $\lim_{t \to 0}t^\la \, x$ exists if $\ell_i \geq 0$ for all $i
\in H$.  Hence a suitable lift $\la$ can be found if 
\begin{equation}
\label{y}
-\al_i^{\vee}=\sum_{j=1}^{N} \ell_j \be_j
\end{equation}
for some $\ell_j \in \Q$ with $\ell_j \geq 0$ whenever $j \in H$.

Since $i \not\in I$, there exists $j \not\in H$ with $\al_i \! \cdot\! \be_j
\neq 0$, indeed $\al_i \!\cdot\! \be_j > 0$ as $\be_j \in \La^+$.  Denote
by $P_i: \La_\Q \to \La_\Q $ the projection onto the hyperplane annihilated by
$\al_i$, that is,
$$P_i(\be_j) = \be_j - \frac{\al_i \!\cdot\! \be_j}{2} \, \al_i^\vee.$$
Rearranging yields
$$-\al_i^\vee = \frac{2}{\al_i \!\cdot\! \be_j}\, P_i(\be_j) -
\frac{2}{\al_i \!\cdot\! \be_j}\, \be_j.$$
By \re{z}, $P_i(\be_j)$ is in the support of $\Si$ and hence
can be expressed as a linear combination of the $\be_i$ with
nonnegative coefficients.  This establishes \re{y}, completing the proof.
\fp

\bit{The wonderful compactification of any semisimple group}

If $G$ is semisimple and the stacky fan $\Si$ comprises only one cone,
the positive Weyl chamber $\La_\Q^+$, equipped with the minimal
cocharacters $\be_i$ along its rays, then we obtain an orbifold whose
coarse moduli space is the compactification proposed by Springer
\cite{springer}.

\bs{Proposition}\label{aa} If $G$ is semisimple and $\Si$ is as
above, then the coarse moduli space $M_G(\Si)$ of $\MM_G(\Si)$ is the
normalization of $M_{G_{\ad}}(\Si) = \overline{G}_{\ad}$ in the function
field of $G$. \es

\bs{Lemma} Suppose given a Cartesian diagram of noetherian integral schemes
\begin{center}
\begin{tikzpicture}
\matrix (n)[matrix of math nodes,  row sep=2em, column sep=1.5em, text height=1.5ex, text depth=0.25ex]
{U & X \\ V & Y, \\};
\path[right hook->] 
(n-1-1) edge (n-1-2);
\path[->]
(n-1-1)  edge (n-2-1);
\path[right hook->] 
(n-2-1) edge  (n-2-2);
\path[->] 
(n-1-2) edge node[auto] {$f$} (n-2-2);
\end{tikzpicture}
\end{center}
where $X$ is normal, $V$ is nonempty and open in $Y$, and $f$ is
finite.  Then $X$ is the normalization of $Y$ in the function field of
$U$. \es

\pf. Let $W$ be the normalization of $Y$ in the function field of $U$.
Since $X$ is normal, $f$ factors through $\tilde{f}:X\to W$, which is
finite because $f$ is.  Since $W$ also contains $U$ as an open
subscheme, $\tilde{f}$ is a finite birational morphism to a normal
scheme, hence an isomorphism by Zariski's main theorem.  \fp

\pf\ of {\rm \re{aa}}.  Take $U = G$, $V=G_{\ad}$,
$X=M_{G}(\La_\Q^+)$, $Y=M_{G_{\ad}}(\La_\Q^+)$ in the lemma.  As a
geometric invariant theory quotient of a normal
variety, $X$ is itself normal.  Since the morphism $\A_\be \to \A$
is finite in this case, so is
$S_{G,\be} \to \SG$ and hence $S_{G,\be} \modmod \G_\be \to S_G \modmod
\G_\be$, that is, $M_G(\La_Q^+) \to
M_{G_{\ad}}(\La_\Q^+)$.  \fp

\bit{The coarse moduli space as a spherical variety}

For any reductive $G$ and simplicial stacky fan $\Si$, the coarse
moduli space $M(\Si)$ may be described as a spherical variety as
follows.  See Timashev \cite{timashev} or Pezzini \cite{pezzini} for
background on spherical varieties.

\bs{Proposition}  
\label{ii}
The coarse moduli space $M(\Si)$ of ${\red \MM_G}(\Si)$ is
the toroidal spherical embedding corresponding to the uncolored fan
$w\Si$, where $w$ is the longest element of $W$.  
\es

Hence any $G \times G$-equivariant toroidal compactification of $G$
with finite quotient singularities is the coarse moduli space of some
${\red \MM_G}(\Sigma)$; for such compactifications correspond to simplicial
fans whose support is the negative Weyl chamber $-u\La_\Q^+$.

\pf.  As the coarse moduli space of a smooth Deligne-Mumford stack,
$M(\Si)$ is normal.  It is also a scheme: indeed, it is covered by the
open sets $M(\si)$, which are schemes since choosing a normal fan
$\Si_\si$ containing $\si$ realizes $M(\si)$ as an open
subset of $M(\Si_\si)$, which is quasiprojective by \re{v}.  Since
$M(\Si)$ contains $G$ as a dense open subset, it also has a dense
orbit for the Borel subgroup $B \times B_- \subset G \times G$, where
$B_-$ is the Borel opposite to $B$.  Hence it is spherical.

To prove it is toroidal, it suffices to show that every $G\times G$-orbit
nontrivially intersects the {\it big cell}, that is, the complement of the
effective divisors preserved by $B \times B_-$ but not $G \times G$.
By \re{u}, it suffices to prove the corresponding statement for the $G
\times G \times \G_\be$-action on $S^0_{G,\be}$.  For $S^0_G$ this is
proved by Vinberg \cite[Prop.\ 14]{v}, and the case
of $S^0_{G,\be}$ immediately follows, since its big cell is the fibered
product of $\A_\be$ with the big cell of $S^0_G$.  

A toroidal spherical embedding is determined by its fan, which is
uncolored and whose support lies in the valuation cone of its dense orbit.
In the present case {\red the dense orbit} is $G = (G \times G)/G$, so the valuation
cone is $-\La_\Q^+$: see Timashev \cite[24.9]{timashev}.  Since by \re{cc} the
closure in $M(\Si)$ of the maximal torus $T$ is the toric variety
$X(W\Si)$, it follows \cite[29.7]{timashev} that the desired fan is
the part of $W \Si$ lying inside the Weyl chamber $-\La_\Q^+$, which
is $w\Si$. \fp

\bit{Relationship with the Losev-Manin space 
\boldmath $\overline{M}_{0,\{1,1,\epsilon,\dots,\epsilon\} }$}

For the root system of type $A_r$, a moduli problem represented by the
toric variety $X(W\La_\Q^+)$ of the fan of Weyl chambers has been
described by Losev-Manin \cite{lm}.  It is best stated in terms of the
{\it weighted pointed stable curves\/} introduced by Hassett
\cite{has}.  These are prestable curves with positive rational weights assigned
to each marked point, where weighted marked points are allowed to
coincide with each other (but not with the nodes) provided that the
sum of their weights does not exceed unity.

\bs{Theorem (Losev-Manin)} The toric variety associated with
the $A_{r}$ root system is the fine moduli space
$\overline{M}_{0,\{1,1,\epsilon,\dots,\epsilon\} }$ of genus
zero weighted pointed stable curves with $r+3$ marked points, two having weight $1$,
and the rest having small weight $\epsilon\ll 1$.  \es

From our point of view these toric varieties arise as the closure of
the maximal torus of $G = PGL_{r+1}$ in $\MM_G(\La_\Q^+)$.  The functor
from Losev-Manin's moduli problem to ours may be described as follows.
Let $C\to S$ be a family of weighted pointed stable curves of the type
described above.  Regard the marked points as sections $p_+, p_-, a_0,
\dots, a_r: S\to C$.  Then $C\to S$ equipped with the sections $p_+$
and $p_-$ is a chain in the sense of \re{bb}.  Let $E\to C$ be the
vector bundle $$E\deq \O(a_0)\oplus\dots\oplus\O(a_r).$$
By \re{i} there is a canonical $\G_m$-action on $C$ lifting to $E$ so that
the action on $p_+^*E$ is trivial.  The action on $p_-^*E$ is then
uniformly of weight $-1$, so the associated $PGL_{r+1}$-bundle is acted on
trivially at $p_\pm$, hence constitutes a framed bundle chain.

Let $\Si$ be the stacky fan whose single cone is the positive Weyl
chamber $\La_\Q^+$ equipped with the fundamental coweights of
$PGL_{r+1}$.  These fundamental coweights $\varpi_i^{\vee}$ are the
minimal cocharacters on the rays of $\La_\Q^+$ and are given by
$\la^{\varpi_i^\vee}(t) = 
\mbox{diag }(t,\dots,t,1,\dots,1)$, 
where $t$ appears $i$ times.  On the other hand, given a framed bundle
of the type described above over a standard chain $C_n$, it is
easily seen that the cocharacter at a node is 
$\mbox{diag }(t^{\ell_0},\dots, t^{\ell_r})$ 
where $\ell_j = -1$ if $a_j$ lies on a component of the chain between
that node and $p_+$, and $\ell_j = 0$ otherwise.  Hence the
cocharacters appearing at the nodes will be a subsequence of
$w(\varpi_1^\vee), \dots, w(\varpi_r^\vee)$ for some fixed permutation
$w \in W$.  The framed bundle chain is therefore $\Si$-stable in the sense
of \re{m}.

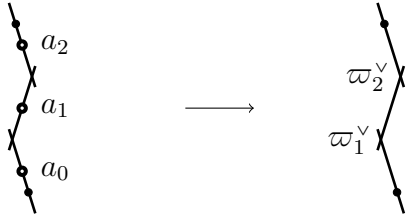
\begin{center}
\begin{tikzpicture}[scale=.7]

\newcommand{\tripleright}[3]
{
\draw[line width=1pt](0,0)--(-.5,1.6);
\draw[line width=1pt](-.5,1.2)--(0,2.8);
\draw[line width=1pt](0,2.4)--(-.5,4);
\filldraw[white] (-.25,2) circle (2pt);
\draw[line width=1.5pt] (-.25,2) circle (2pt);
\draw (-.15, 2) node[anchor=  west] {$#1$};
\filldraw[white] (-.25,.8) circle (2pt);
\draw[line width=1.5pt] (-.25,.8) circle (2pt);
\draw (-.15, .8) node[anchor=  west] {$#2$};
\filldraw[white] (-.25,3.2) circle (2pt);
\draw[line width=1.5pt] (-.25,3.2) circle (2pt);
\draw (-.15, 3.2) node[anchor=  west] {$#3$};
\filldraw (-0.13,.4) circle (2pt);
\filldraw (-0.37,3.6) circle (2pt);
\draw (-.03, .4) node[anchor=  west] {};
\draw (-.27, 3.6) node[anchor=  west] {};
}

\begin{scope}[xshift=-3.5cm, yshift=-2cm]
\tripleright{a_1}{a_0}{a_2}
\end{scope}

\begin{scope}[scale=1.3]
\draw[->] (-.5cm,0) -- (.5cm,0);
\end{scope}

\begin{scope}[xshift=3.5cm, yshift=-2cm]

\draw[line width=1pt](0,0)--(-.5,1.6);
\draw[line width=1pt](-.5,1.2)--(0,2.8);
\draw[line width=1pt](0,2.4)--(-.5,4);

\filldraw (-0.13,.4) circle (2pt);
\filldraw (-0.37,3.6) circle (2pt);
\draw (-.03, .4) node[anchor=  west] {};
\draw (-.27, 3.6) node[anchor=  west] {};

\draw (-1.7, 1.4) node[anchor=  west] {$\varpi_1^{\vee}$};
\draw (-1.35, 2.6) node[anchor=  west] {$\varpi_2^{\vee}$};

\end{scope}

\end{tikzpicture}
\end{center}

This construction provides a functor from weighted pointed stable
curves to framed bundle chains, yielding a morphism
$\overline{M}_{0,\{1,1,\epsilon,\dots,\epsilon\} } \to \MM_G(\Si)$.
Since the framings respect the splitting into line bundles, this
actually factors through $\MM_T(W\Si)$.  Indeed, by arguing as in the
proof of \re{cc} one shows that
$\overline{M}_{0,\{1,1,\epsilon,\dots,\epsilon\} } \to
\MM_T(W\Si)$ is an isomorphism.

\begin{figure}[b]

\begin{center}
\rule{5cm}{.1mm}
\end{center}

\begin{center}
\begin{tikzpicture}[scale=.7]

\newcommand{\bottomdoubleright}[3]
{\draw[line width=1pt](0,0)-- (.5,2.4); 
\draw[line width=1pt](.5,1.6)-- (0,4); 
\filldraw (.078,.4) circle (2pt);
\filldraw (.078,3.6) circle (2pt);
\draw (0.178, .4) node[anchor=   west] {};
\draw (0.178, 3.6) node[anchor=  west] {};
\filldraw[white] (.2,3) circle (2pt);
\draw[line width=1.5pt] (.2,3) circle (2pt);
\draw (0.3, 3) node[anchor=  west] {\small \it #1};
\filldraw[white] (.165,.8) circle (2pt);
\draw[line width=1.5pt] (.165,.8) circle (2pt);
\draw (0.265, .8) node[anchor=  west] {\small \it #2};
\filldraw[white] (.27,1.3) circle (2pt);
\draw[line width=1.5pt] (.27,1.3) circle (2pt);
\draw (0.37, 1.3) node[anchor=  west] {\small \it #3};
}

\newcommand{\bottomdoubleleft}[3]
{\draw[line width=1pt](0,0)-- (.5,2.4); 
\draw[line width=1pt](.5,1.6)-- (0,4); 
\filldraw (.078,.4) circle (2pt);
\filldraw (.078,3.6) circle (2pt);
\draw (-.022, .4) node[anchor=  east] {};
\draw (-0.022, 3.6) node[anchor=  east] {};
\filldraw[white] (.2,3) circle (2pt);
\draw[line width=1.5pt] (.2,3) circle (2pt);
\draw (0.1, 3) node[anchor=  east] {\small \it #1};
\filldraw[white] (.165,.8) circle (2pt);
\draw[line width=1.5pt] (.165,.8) circle (2pt);
\draw (0.065, .8) node[anchor=  east] {\small \it #2};
\filldraw[white] (.27,1.3) circle (2pt);
\draw[line width=1.5pt] (.27,1.3) circle (2pt);
\draw (0.17, 1.3) node[anchor=  east] {\small \it #3};
}

\newcommand{\topdoubleright}[3]
{
\draw[line width=1pt](0,0)-- (-.5, 2.4); 
\draw[line width=1pt](-.5,1.6)-- (0,4); 

\filldraw (-0.078,.4) circle (2pt);
\filldraw (-0.078,3.6) circle (2pt);
\draw (.022, .4) node[anchor=  west] {};
\draw (0.022, 3.6) node[anchor=  west] {};
\filldraw[white] (-.175,3.1) circle (2pt);
\draw[line width=1.5pt] (-.175,3.1) circle (2pt);
\draw (-0.075, 3.1) node[anchor=  west] {\small \it #1};
\filldraw[white] (-.29,2.6) circle (2pt);
\draw[line width=1.5pt] (-.29,2.6) circle (2pt);
\draw (-.19, 2.6) node[anchor=  west] {\small \it #2};
\filldraw[white] (-.27,1.3) circle (2pt);
\draw[line width=1.5pt] (-.27,1.3) circle (2pt);
\draw (-.17, 1.4) node[anchor=  west] {\small \it #3};
}

\newcommand{\topdoubleleft}[3]
{
\draw[line width=1pt](0,0)-- (-.5, 2.4); 
\draw[line width=1pt](-.5,1.6)-- (0,4); 

\filldraw (-0.078,.4) circle (2pt);
\filldraw (-0.078,3.6) circle (2pt);
\draw (-.178, .4) node[anchor=  east] {};
\draw (-.178, 3.6) node[anchor=  east] {};
\filldraw[white] (-.175,3.1) circle (2pt);
\draw[line width=1.5pt] (-.175,3.1) circle (2pt);
\draw (-.275, 3.1) node[anchor=  east] {\small \it #1};
\filldraw[white] (-.29,2.6) circle (2pt);
\draw[line width=1.5pt] (-.29,2.6) circle (2pt);
\draw (-.39, 2.6) node[anchor=  east] {\small \it #2};
\filldraw[white] (-.27,1.3) circle (2pt);
\draw[line width=1.5pt] (-.27,1.3) circle (2pt);
\draw (-.37, 1.4) node[anchor=  east] {\small \it #3};
}

\newcommand{\tripleleft}[3]
{
\draw[line width=1pt](0,0)--(-.5,1.6);
\draw[line width=1pt](-.5,1.2)--(0,2.8);
\draw[line width=1pt](0,2.4)--(-.5,4);
\filldraw[white] (-.25,2) circle (2pt);
\draw[line width=1.5pt] (-.25,2) circle (2pt);
\draw (-.35, 2) node[anchor=  east] {\small \it #1};
\filldraw[white] (-.25,.8) circle (2pt);
\draw[line width=1.5pt] (-.25,.8) circle (2pt);
\draw (-.35, .8) node[anchor=  east] {\small \it #2};
\filldraw[white] (-.25,3.2) circle (2pt);
\draw[line width=1.5pt] (-.25,3.2) circle (2pt);
\draw (-.35, 3.2) node[anchor=  east] {\small \it #3};
\filldraw (-0.13,.4) circle (2pt);
\filldraw (-0.37,3.6) circle (2pt);
\draw (-.23, .4) node[anchor=  east] {};
\draw (-.47, 3.6) node[anchor=  east] {};

}

\newcommand{\tripleright}[3]
{
\draw[line width=1pt](0,0)--(-.5,1.6);
\draw[line width=1pt](-.5,1.2)--(0,2.8);
\draw[line width=1pt](0,2.4)--(-.5,4);
\filldraw[white] (-.25,2) circle (2pt);
\draw[line width=1.5pt] (-.25,2) circle (2pt);
\draw (-.15, 2) node[anchor=  west] {\small \it #1};
\filldraw[white] (-.25,.8) circle (2pt);
\draw[line width=1.5pt] (-.25,.8) circle (2pt);
\draw (-.15, .8) node[anchor=  west] {\small \it #2};
\filldraw[white] (-.25,3.2) circle (2pt);
\draw[line width=1.5pt] (-.25,3.2) circle (2pt);
\draw (-.15, 3.2) node[anchor=  west] {\small \it #3};
\filldraw (-0.13,.4) circle (2pt);
\filldraw (-0.37,3.6) circle (2pt);
\draw (-.03, .4) node[anchor=  west] {};
\draw (-.27, 3.6) node[anchor=  west] {};

}

\filldraw[gray!15!white](75:5)--(105:5)--(195:5)--(225:5)--(315:5)--(345:5)--cycle;
\draw[gray, line width=3pt](75:5)--(105:5)--(195:5)--(225:5)--(315:5)--(345:5)--cycle;

\begin{scope}[xshift=1cm, yshift= 0]

\draw[line width=1pt](0,2)--(0,-2);
\filldraw[white] (0,.8) circle (2pt);
\draw[line width=1.5pt] (0,.8) circle (2pt);
\draw (0.1, .8) node[anchor=  west] {\small \it a};
\filldraw[white] (0,0) circle (2pt);
\draw[line width=1.5pt] (0,0) circle (2pt);
\draw (0.1, 0) node[anchor=  west] {\small \it b};
\filldraw[white] (0,-.8) circle (2pt);
\draw[line width=1.5pt] (0,-.8) circle (2pt);
\draw (0.1, -.8) node[anchor=  west] {\small \it c};
\filldraw (0,-1.6) circle (2pt);
\draw (0.1, -1.6) node[anchor=  west] {};
\filldraw (0,1.6) circle (2pt);
\draw (0.1, 1.6) node[anchor=  west] {};
\end{scope}

\begin{scope}[xshift=6cm, yshift=-5.5cm]
\bottomdoubleright{a}{c}{b}
\end{scope}

\begin{scope}
[xshift=-6.4cm, yshift=-5.5cm]
\bottomdoubleleft{c}{b}{a}
\end{scope}

\begin{scope}
[xshift=-.1cm, yshift=6.5cm]
\bottomdoubleright{b}{c}{a}
\end{scope}

\begin{scope}[xshift=5cm, yshift=.5cm]
\topdoubleright{a}{b}{c}
\end{scope}

\begin{scope}[xshift=-5cm,yshift=.5cm]
\topdoubleleft{b}{c}{a}
\end{scope}

\begin{scope}[yshift=-8.5cm]
\topdoubleleft{a}{c}{b}
\end{scope}

\begin{scope}[xshift=3.5cm,yshift=5.5cm]
\tripleright{a}{c}{b}
\end{scope}

\begin{scope}[xshift=-3.5cm,yshift=5.5cm]
\tripleleft{c}{a}{b}
\end{scope}

\begin{scope}[xshift=8.5cm,yshift=-2cm]
\tripleright{b}{c}{a}
\end{scope}

\begin{scope}[xshift=-8.5cm,yshift=-2cm]
\tripleleft{b}{a}{c}
\end{scope}

\begin{scope}[xshift=-3.5cm,yshift=-9cm]
\tripleright{a}{b}{c}
\end{scope}

\begin{scope}[xshift=3.5cm,yshift=-9cm]
\tripleleft{c}{b}{a}
\end{scope}

\draw[very thick,gray,->](4,-7)..controls (5,-7) and (5,-6) .. (315:5.1);
\draw[very thick,gray,->](-4,-7)..controls (-5,-7) and (-5,-6) .. (225:5.1);

\draw[very thick,gray,->](.5,-7)..controls (1.5,-7) and (1.5,-6) .. (0,-3.7);

\draw[very thick,gray,->](6,-3.6)..controls (5.5,-3.8) and (5,-4) .. (330:5);
\draw[very thick,gray,->](-5.7,-3.6)..controls (-5.5,-3.6) and (-5,-4) .. (210:5);

\draw[very thick,gray,->](.7,0)..controls (.3,-.6) and (-1.8,-2) .. (-1,0);

\draw[very thick,gray,->](8,0)..controls (7,0) and (6,-.5) .. (345:5.1);
\draw[very thick,gray,->](-8.5,0)..controls (-7,0) and (-6,-.5) .. (194:5.1);

\draw[very thick,gray,->](4.3,2.5)..controls (4,2.5) and (3.5,2.3) .. (30:3.65);
\draw[very thick,gray,->](-5.1,2.5)..controls (-4,2.5) and (-3.5,2.3) .. (150:3.65);

\draw[very thick,gray,->](2.7,7.5)..controls (2.2,7.5) and (1.5,6) .. (75:5.1);
\draw[very thick,gray,->](-3,7.5)..controls (-2.5,7.5) and (-1.5,6) .. (105:5.1);

\draw[very thick,gray,->](0,8.5)..controls (-1,8.5) and (-.5,7) .. (90:5);
\vspace{1cm}

\
\end{tikzpicture}
\end{center}

\noindent \emph{\small The moment map for the toric variety
  $\overline{\mathcal{M}}_{0,\{1,1,\epsilon,\epsilon,\epsilon\}}
\cong\MM_T(W\La_\Q^+)$,
with image the permutohedron of the symmetric group $S_3$. 
The points $p_\pm =$ \begin{tikzpicture}[scale=.7]
\filldraw (.2,3) circle (2pt);
\end{tikzpicture} have weight 1; the points $a,b,c=$ \begin{tikzpicture}[scale=.7]
\draw[line width=1.5pt] (0,0) circle (2pt); \end{tikzpicture} have weight $\epsilon$.
}

\end{figure}

The Losev-Manin picture was partially extended to other classical
groups by Batyrev and Blume \cite{batblum}.  The correspondence works
best for the $B_r$ root system, where the toric variety can be
interpreted as a moduli space of weighted pointed stable curves with
an involution.  From our point of view this is because the classical
group corresponding to the $B_r$ root systems, namely $SO_{2r+1}$, has
trivial center and hence can be embedded in $PGL_{2r+1}$, indeed as the
fixed-point locus of an involution.  This is not the case with
$SO_{2n}$ and $Sp_n$, leading to the complications observed in those
cases by Batyrev-Blume.

\bit{Relationship with the Kausz space \boldmath $KGL_r$}

Kausz \cite{kausz} has introduced an equivariant compactification
$KGL_r$ of $GL_r$.  It represents the moduli problem of
``bf-morphisms,'' also introduced by Kausz.  There is presumably a
faithful functor from bf-morphisms to bundle chains.  However, since
Kausz's moduli problem is quite technical, we prefer to study $KGL_r$ via its
alternative description by Kausz as an iterated blow-up of
projective space, parallel to Vainsencher's
description \cite{vainsencher} of the wonderful compactification
$\overline{PGL}_r$.

Recall first Vainsencher's description.  Let $X_0 = \P^{r^2-1}$ be the
projectivization of the vector space of $r \times r$ matrices, and let
$A_i$ be the locus of matrices of rank $i$.  Recursively define $X_i
= \Bl(X_{i-1},\overline{A}_i)$.  Vainsencher proves that each blow-up
locus, and hence each $X_i$, is smooth; and that moreover
$X_{r-1} \cong \overline{PGL}_r$, the wonderful compactification.
Furthermore, there is a recursive structure: each $A_i$ is plainly a
$PGL_i$-bundle over a product of Grassmannians (the projection being given
by kernel and image), and its closure $\overline{A}_i$ in
$X_{i-1}$ is the associated $\overline{PGL}_i$-bundle.

Kausz's description is similar but one dimension greater.  Let $Y_0
= \P^{r^2}$ be regarded as compactifying the vector space $\A^{r^2}$
of all $r \times r$ matrices, and refer to elements of the hyperplane
$\P^{r^2-1} = \P^{r^2} \sans \A^{r^2}$ as {\it matrices at infinity}.
Let $B_i \subset \A^{r^2}$ be the locus of matrices of rank $i$, and
let $C_i \subset \P^{r^2-1}$ be the locus of matrices at infinity of
rank $i$.  Recursively define $Y_i
= \Bl(Y_{i-1},\overline{B}_{i-1} \cup \overline{C}_i)$.  Kausz proves
that each blow-up locus is a disjoint union of two smooth loci, so
that $Y_i$ is also smooth; and that moreover $Y_{r-1} \cong KGL_r$,
the space representing his moduli problem.  Furthermore, there is
a recursive structure: each $B_i$ is a $GL_i$-bundle over a product of
Grassmannians, and its closure in $Y_i$ is the associated
$KGL_i$-bundle; while each $C_i$ is a $PGL_i$-bundle over a product of
Grassmannians, and its closure in $Y_{i-1}$ is the associated
$\overline{PGL}_i$-bundle.


For example, $KGL_2$ is the blow-up of $\P^4$, with projective
coordinates $[a,b,c,d,e]$, at the point $[0,0,0,0,1]$ and the quadric
surface $ad-bc=0=e$.  

In fact, the Vainsencher and Kausz constructions are intimately related.  
For assigning $A \mapsto \P(1 \oplus A)$ gives an embedding $GL_r \to
PGL_{r+1}$ restricting to an isomorphism on the standard maximal tori.
It extends to a linear inclusion $Y_0 \subset X_0$ of the projective
space compactifications (with $r+1$ substituted for $r$ in the case of
$X_0$).  It is obvious that $A_i \cap Y_0 = B_{i-1} \cup C_i$ in
$X_0$.  From a careful scrutiny of the recursive structure, one should
be able to see further that $\overline{A}_i \cap Y_{i-1}
= \overline{B}_{i-1} \cap \overline{C}_i$ in $X_{i-1}$, even
scheme-theoretically, and hence, blowing up both sides, that
$Y_i \subset X_i$ as the closure of $GL_r$.  Consequently, $KGL_r$ is
nothing but the closure of $GL_r$ in the wonderful compactification
$\overline{PGL}_{r+1}$.  In any case, the latter statement is proved,
using the classification of toroidal group compactifications, in the
thesis of Huruguen \cite[3.1.16]{huruguen-thesis}.

Kausz's space may be realized as a moduli space of bundle chains.
This is assured by \re{ii}, as it is a smooth toroidal equivariant
compactification of $GL_r$.  The corresponding stacky fan $\Si$ may be
described as follows.  As the maximal tori of $GL_r$ and $PGL_{r+1}$
are identified, the cocharacter lattices of the two groups may be
identified as well.  Let $\Si$ be that part of the fan of Weyl
chambers of $PGL_{r+1}$ lying in the positive Weyl chamber of $GL_r$,
equipped with the minimal cocharacters $\be_i$ along its rays.  For
example, when $r=2$, the fan $\Si$ consists of three contiguous
$60^\circ$ sectors in the plane.

\bs{Proposition} With notation as above, $KGL_r \cong {\red \MM_{GL_r}}(\Si)$. \es

\pf.  Since the Weyl groups satisfy $W_{GL_r} \subset W_{PGL_{r+1}}$,
the fan $W_{GL_r} \Si$ is precisely the fan of all Weyl chambers of
$PGL_r$.  Extension of structure group for bundle chains then
determines a morphism ${\red \MM_{GL_r}} \to \overline{PGL}_{r+1}$.  This is an
embedding: it separates points since a (framed) vector bundle $E$ may
be recovered from the (framed) projective bundle $\P(E \oplus \co)$,
and it separates tangent vectors as is easily seen from the
description of the Zariski tangent spaces in \re{q}.  Hence ${\red \MM_{GL_r}}$
is a subscheme of $\overline{PGL}_{r+1}$ containing $GL_r$ as a dense
open subset.  It therefore coincides with $KGL_r$, the closure of
$GL_r$ in $\overline{PGL}_{r+1}$.  \fp

\begin{figure}[t]
\begin{center}
\begin{tikzpicture}

\begin{scope}[xshift=-1.8cm, yshift=0cm, scale=1.3]

\draw[very thick] (-4,2) -- (3,2);
\draw[very thick] (-4,-2) -- (3,-2);
\draw[very thick] (1.5,-2.6) -- (3.0,0.6);
\draw[very thick] (3.0,-0.6) -- (1.5,2.6);

\draw (-2.5,0) node[anchor = center] {$\O \oplus \O$};
\draw (-2.0,2) node[anchor = south] {$\O(-1,1) \oplus \O(-1,1)$};
\draw (-2.0,-2) node[anchor = north] {$\O(1,-1) \oplus \O(1,-1)$};
\draw (2.9,0) node[anchor = west] {$\leftarrow \, \O(1,-1,0) \oplus \O(0,-1,1)$};
\draw (3.8,2) node[anchor = south] {$\swarrow \, \O(-1,1,0) \oplus \O (-1,0,1)$};
\draw (3.8,-2) node[anchor = north] {$\nwarrow \, \O(1,0,-1) \oplus \O (0,1,-1)$};
\draw (2.3,1.2) node[anchor = west] {$\O(0,0) \oplus \O(-1,1)$};
\draw (2.3,-1.2) node[anchor = west] {$\O(1,-1) \oplus \O(0,0)$};

\end{scope}

\end{tikzpicture}
\end{center}

\noindent \emph{\small 
The orbits of the $GL_2 \times GL_2$-action on $KGL_2$, labeled with
the isomorphism classes of their corresponding bundle chains, using the multidegree notation indicated in the text.  }

\begin{center}
\rule{5cm}{.1mm}
\end{center}

\end{figure}

Let us spell out concretely which bundle chains are stable in the moduli
problems represented by $\overline{PGL}_r$ and $KGL_r$.  All the
relevant bundles have trivial weights at $p_\pm$.  So to specify a
line bundle over the standard chain $C_n$, rather than use the notation
$\O(b_0 \st b_1 \st \cdots \st b_{n+1})$ of \re{c}(b), it is now more
convenient to use the multidegree notation $\O(d_0, \dots, d_n)$,
where $d_i = b_i - b_{i+1}$ is the degree on the $i$th component of
the chain.  

With this notation, the bundle chains over the standard chain $C_n$
which are stable in the moduli problem represented by
$\overline{PGL}_r$ are those of the form
$\P(\O(v_1) \oplus \cdots \oplus \O(v_r))$, where each $v_j$ is a
standard basis vector in $\Z^{n+1}$, each of which must appear at
least once.  For example, the only bundle allowed over $C_0 = \P^1$ is
$\P(\O(1) \oplus \cdots \oplus \O(1))$, the trivial bundle.

Likewise, the bundle chains over the standard chain $C_n$ which are
stable in the moduli problem represented by $KGL_r$ are those vector
bundles $E$ such that $\P(\O \oplus E)$ is of the form stated in the
last paragraph.  Equivalently,
$E \cong \O(v_1-v_0) \oplus \cdots \oplus \O(v_r-v_0)$, where each
$v_j$ is a standard basis vector in $\Z^{n+1}$, each of which must
appear at least once.

For example, $KGL_2$, has 8 $GL_2 \times GL_2$-orbits: 1 of
codimension 0, 4 of codimension 1, and 3 of codimension 2.  Their
incidences, and the isomorphism classes of the vector bundles appearing in the
corresponding bundle chains, are shown in the accompanying diagram.

From the fan perspective, the choice of $KGL_r$ as a compactification
of $GL_r$ seems somewhat arbitrary.  There are many other possible
fans covering the Weyl chamber of $GL_r$ leading to equally good
compactifications.  For example, when $r=2$, the fan consisting of two
contiguous $90^\circ$ sectors in the plane leads to the
compactification $[(\overline{SL}_2 \times \P^1)/\mu_2]$.

On the other hand, several of the alternative descriptions of
$\overline{PGL}_r$ given by the second author \cite{ccr} carry over
neatly to $KGL_r$ by tracing through the closures of the relevant
loci, starting from $GL_r \subset PGL_{r+1}$.  For example, let $\G_m$
act by scalar multiplication on an $r$-dimensional vector space $V$
and thus on the Grassmannian $X = \Gr_r(V \oplus V^*)$.  Then
$\overline{PGL}_r$ is the Chow quotient of $X$ by $\G_m$; and
likewise, $KGL_r$ is the Chow quotient of $\P^1 \times X$ by $\G_m$.
Or again, if $\ev: \overline{M}_{0,2}(X,r) \to X^2$ denotes evaluation
on the moduli space of stable maps, then $\overline{PGL}_r
= \ev^{-1}(V, V^*)$; and likewise, $KGL_r$ is a
fiber of evaluation on the so-called ``graph space''.  That is to say, if
$\Ev: \overline{M}_{0,2}(\P^1 \times X, 1 \times r) \to (\P^1 \times
X)^2$ again denotes evaluation, then $KGL_r = \Ev^{-1}({[1,0] \times
V}, \, {[0,1] \times V^*})$.  From this perspective, the
choice of $KGL_r$ seems less arbitrary.

\end{document}